\documentclass[12pt]{amsart}
\usepackage{amssymb, eucal, amsfonts, amsmath, xypic,latexsym}

\def\gr{\mathrm{gr}}
\def\k{{\Bbbk}}

\def\g{{\mathfrak g}}

\def\h{{\mathfrak h}}
\def\sl{{\mathfrak sl}}

\def\z{{\mathfrak z}}
\def\m{{\mathfrak m}}
\def\n{{\mathfrak n}}

\def\Z{{\mathbb Z}}

\def\End{\mathop{\fam0 End}}

\def\Lie{\mathop{\fam0 Lie}}
\def\Ad{\mathrm{Ad}\,}

\def\Mat{\mathop{\fam0 Mat}\nolimits}
\def\sl{\mathop{\mathfrak{sl}}\nolimits}

\def\ad{\mathrm{ad\,}}

\def\la{\langle}
\def\ra{\rangle}

\theoremstyle{plain}
\newtheorem{theorem}{Theorem}[section]
\newtheorem{corollary}{Corollary}[section]
\newtheorem{prop}{Proposition}[section]
\newtheorem{lemma}{Lemma}[section]

\theoremstyle{definition}

\theoremstyle{remark}

\newtheorem{rem}{Remark}[section]

\def\subtitle#1. {{\medskip\bf#1\par\nobreak\smallskip}}
\def\proclaim#1. {\medbreak\bgroup\noindent\bf#1. \it}

\def\endproclaim{\egroup
\ifdim\lastskip<\medskipamount\removelastskip\medskip\fi}
\newcount
=0
\def\citedef#1 {\advance\citation by1
  \expandafter\edef\csname#1\endcsname{{\the\citation}}
  \checkendcitedef}
\def\checkendcitedef#1{\ifx#1\endcitedef\else\citedef#1\fi}
\def\cite#1{\csname#1\endcsname}
\citedef BV1 BV2 BK Bo Bou BrK DDeC DeSK  Du Eis GG G Ja Ja1 J84
Kr P95 P02  P03 P07  Sk SpSt
\endcitedef
\newtoks\nextauth
\newif\iffirstauth
\def\checkendauth#1{\ifx\endauth#1
        \iffirstauth\the\nextauth
        \else{} and \the\nextauth\fi,
    \else\iffirstauth\the\nextauth\firstauthfalse
        \else, \the\nextauth\fi
        \expandafter\auth\expandafter#1\fi}
\def\auth#1 #2 {\nextauth={#1 #2}\checkendauth}
\newif\ifinbook
\newif\ifbookref
\def\nextref#1 {\bookreffalse\inbookfalse
    \bibitem[\cite{#1}]{}
    \firstauthtrue
    \ignorespaces}
\def\paper#1{{\it#1,}}
\def\In#1{\inbooktrue In #1,}
\def\book#1{\bookreftrue{\it#1,}}
\def\journal#1{#1\ifinbook,\fi}
\def\bookseries#1{#1,}
\def\Vol#1{\ifbookref Vol. #1,\else\ifinbook Vol. #1,\else{\bf#1}\fi\fi
    \space\ignorespaces}

\def\publisher#1{#1,}
\def\Year#1{\ifbookref #1.\else\ifinbook #1,\else(#1)\fi\fi
    \space\ignorespaces}
\def\Pages#1{\ifinbook pp. #1.\else #1.\fi}
\begin{document}
\title{Primitive ideals, non-restricted representations and finite $W$-algebras}
\author{Alexander Premet}
\thanks{\nonumber{\it Mathematics Subject Classification} (2000 {\it revision})).
Primary 17B35. Secondary 17B63, 17B81. {\it \mbox{\,\quad}Key
words and phrases}: primitive ideal, modular representation,
finite $W$-algebra}
\address{School of Mathematics, University of Manchester, Oxford Road,
M13 9PL, UK} \email{sashap@maths.man.ac.uk}
\begin{abstract}
\noindent Let $G$ be a simple algebraic group over $\mathbb C$ and
$\g=\text{Lie}\,G$. Let $(e,h,f)$ be an $\sl_2$-triple in $\g$ and
$(\,\cdot\,,\,\cdot\,)$ the $G$-invariant bilinear form on $\g$
such that $(e,f)=1$. Let $\chi\in\g^*$ be such that
$\chi(x)=(e,x)$ for all $x\in\g$ and let $H_\chi$ denote the
enveloping algebra of the Slodowy slice $e+{\rm Ker}\,\ad f$. Let
$\mathcal I$ be a primitive ideal of the universal enveloping
algebra $U(\g)$ whose associated variety is the closure of the
coadjoint orbit $\mathcal O_\chi$. We prove in this note that if
$\mathcal I$ has rational infinitesimal character, then there is a
finite dimensional irreducible $H_\chi$-module $V$ such that
${\mathcal I}={\rm Ann}_{U(\g)}\big(Q_\chi\otimes_{H_\chi}V\big)$,
where $Q_\chi$ is the generalised Gelfand--Graev $\g$-module
associated with the triple $(e,h,f)$. In conjunction with
well-known results of Barbasch and Vogan this implies that all
finite $W$-algebras possess finite dimensional irreducible
representations.
\end{abstract}

\maketitle





\section{\bf Introduction}
\subsection{}\label{1.1}
Let $G$ be a simple, simply connected algebraic group over $\mathbb
C$. Let $\g=\Lie(G)$ and let $(e,h,f)$ be an $\sl_2$-triple in $\g$.
Let $(\,\cdot\,,\,\cdot\,)$ be the $G$-invariant bilinear form on
$\g$ with $(e,f)=1$ and define $\chi=\chi_e\in\g^*$ by setting
$\chi(x)=(e,x)$ for all $x\in\g$. Denote by ${\mathcal O}_\chi$ the
coadjoint $G$-orbit of $\chi$.

Let ${\mathcal S}_e=e+\text{Ker}\,\ad\,f$ be the Slodowy slice at
$e$ through the adjoint orbit of $e$ and let $H_\chi$ be the
enveloping algebra of ${\mathcal S}_e$; see [\cite{P02},
\cite{GG}, \cite{BrK}]. Recall that
$H_\chi={\End}_{\g}\,(Q_\chi)^{\text{op}}$ where $Q_\chi$ is the
generalised Gelfand--Graev module for $U(\g)$ associated with the
triple $(e,h,f)$. The $\g$-module $Q_\chi$ is induced from a
one-dimensional module ${\mathbb C}_\chi$ over of a nilpotent
subalgebra $\m$ of $\g$ such that $\dim \m=\frac{1}{2}\dim
{\mathcal O}_\chi$. The subalgebra $\m$ is $(\ad h)$-stable, all
weights of $\ad h$ on $\m$ are negative, and $\chi$ vanishes on
$[\m,\m]$. The action of $\m$ on ${\mathbb C}_\chi={\mathbb C}
1_\chi$ is given by $x(1_\chi)=\chi(x)1_\chi$ for all $x\in\m$;
see [\cite{P02}, \cite{GG}, \cite{BrK}, \cite{P07}] for more
detail.

Denote by $\z_\chi$ the stabiliser of $\chi$ in $\g$ (it coincides
with the centraliser $\g_e$ of $e$ in $\g$). By the $\sl_2$-theory,
all weights of $\ad h$ on $\z_\chi$ are nonnegative integers. Let
$x_1,\ldots,x_r$ be a basis of $\z_\chi$ such that $[h,x_i]=n_i x_i$
for some $n_i\in\Z_+$. By [\cite{P02}], to each basis vector $x_i$
one can attach an element $\Theta_{x_i}\in H_\chi$ in such a way
that the monomials $\Theta_{x_1}^{i_1}\Theta_{x_2}^{i_2}\cdots
\Theta_{x_r}^{i_r}$ with $(i_1,i_2,\ldots,i_r)\in \Z_+^r$ form a
basis of $H_\chi$ over $\mathbb C$. The monomial
$\Theta_{x_1}^{a_1}\Theta_{x_2}^{a_2}\cdots \Theta_{x_r}^{a_r}$ is
said to have {\it Kazhdan degree} $\sum_{i=1}^r a_i(n_i+2)$. For
$k\ge 0$, we denote by $H_\chi^k$ the $\mathbb C$-span of all
monomials as above of Kazhdan degree $\le k$. By [\cite{P02}, 4.6],
we have  $H_\chi^i\cdot H_\chi^j\subseteq H_\chi^{i+j}$ for all
$i,j\in\Z_+$. Thus $\{H_\chi^k\,|\,k\in\Z_+\}$ is an increasing
filtration of the algebra $H_\chi$. The corresponding graded algebra
$\text{gr}\,H_\chi$ is a polynomial algebra in
$\text{gr}\,\Theta_{x_1},\text{gr}\,\Theta_{x_2},\ldots,\text{gr}\,
\Theta_{x_r}$, and it identifies naturally with the coordinate
algebra ${\mathbb C}[{\mathcal S}_e]$ endowed with its Slodowy
grading. The Poisson bracket on ${\mathbb C}[{\mathcal S}_e]\cong
\gr\,H_\chi$ induced by multiplication in $H_\chi$ is determined by
Gan--Ginzburg in [\cite{GG}].

According to recent results of D'Andrea--De Concini--De
Sole--Heluany--Kac [\cite{DDeC}] and De Sole--Kac [\cite{DeSK}],
the algebra $H_\chi$ is isomorphic to the Zhu algebra of the
vertex $W$-algebra associated with $\g$ and $e$. The latter
algebra is, in turn, isomorphic to the finite $W$-algebra $W^{\rm
fin}(\g,e)$ obtained from the Poisson algebra $\gr\, H_\chi$ by
BRST quantisation; see [\cite{DeSK}] for detail. Thus,
$H_\chi\cong W^{\rm fin}(\g,e)$ as predicted in [\cite{P02},
1.10].

\subsection{}\label{1.2}
Let ${\mathcal C}_\chi$ denote the category of all $\g$-modules on
which $x-\chi(x)$ acts locally nilpotently for every $x\in\m$.
Given a $\g$-module $M$ we set $$\text{Wh}(M):=\,\{m\in
M\,|\,x.m=\chi(x)m\ \ \,\forall\,x\in\m\}.$$ The algebra $H_\chi$
acts on $\text{Wh}(M)$ via a canonical isomorphism $H_\chi\cong
\big(U(\g)/N_\chi\big)^{\ad \m}$ where $N_\chi$ denotes the left
ideal of the universal enveloping $U(\g)$ generated by all
$x-\chi(x)$ with $x\in\m$. By Skryabin's theorem [\cite{Sk}], the
functors $V\rightsquigarrow Q_\chi\otimes_{H_\chi}V$ and
$M\rightsquigarrow \mbox{Wh}(M)$ are mutually inverse equivalences
between the category of all $H_\chi$-modules and the category
${\mathcal C}_\chi$; see also [\cite{GG}, Theorem~6.1]. Skryabin's
equivalence implies that for any irreducible $H_\chi$-module $V$
the annihilator $\text{Ann}_{U(\g)}(Q_\chi\otimes_{H_\chi} V)$ is
a primitive ideal of $U(\g)$. By the Irreducibility Theorem, the
associated variety $\mathcal{VA}(I)$ of any primitive ideal $I$ of
$U(\g)$ is the closure of a nilpotent orbit in $\g^*$; see
[\cite{G}] and references therein.

Given a finitely generated $U(\g)$-module $M$ we denote by
$\text{Dim}(M)$ the Gelfand--Kirillov dimension of $M$. By
[\cite{P07}, Theorem~3.1], for any irreducible $H_\chi$-module $V$
the associated variety of $\text{Ann}_{U(\g)}(Q_\chi\otimes_{H_\chi}
V)$ contains ${\mathcal O}_\chi$, and if $\dim V<\infty,$  then
$$\mathcal{VA}\big(\text{Ann}_{U(\g)}(Q_\chi\otimes_{H_\chi} V)\big)
=\overline{\mathcal O}_\chi\quad\,\text{and}\,\quad
\text{Dim}(Q_\chi\otimes_{H_\chi} V)=\frac{1}{2}\dim
\overline{\mathcal O}_\chi.$$ In particular, if the
$H_\chi$-module $V$ is finite dimensional, then the simple
$U(\g)$-module $Q_\chi\otimes_{H_\chi} V$ is holonomic. It was
conjectured in [\cite{P07}, 3.4] that for any primitive ideal
$\mathcal I$ of $U(\g)$ with $\mathcal{VA}({\mathcal
I})=\overline{\mathcal O}_\chi$ there exists a finite dimensional
irreducible $H_\chi$-module $V$ such that ${\mathcal
I}=\text{Ann}_{U(\g)}(Q_\chi\otimes_{H_\chi} V)$. In [\cite{P07},
5.6], this conjecture was proved under the assumption that $e$
belongs to the minimal nilpotent orbit of $\g$.

Let $\h$ be a Cartan subalgebra of $\g$, and $\Phi$ the root system
of $\g$ relative to $\h$. Let $\Pi=\{\alpha_1,\ldots,\alpha_\ell\}$
be a basis of simple roots in $\Phi$ with the elements in $\Pi$
numbered as in [\cite{Bou}], and $\Phi^+$ the positive system of
$\Phi$ relative to $\Pi$. Let
$\rho=\frac{1}{2}\sum_{\alpha\in\Phi^+}\alpha$ and let $W$ be the
Weyl group of $\Phi$. Given $\lambda\in\h^*$ we let $L(\lambda)$
denote the irreducible $\g$-module with highest weight $\lambda$. As
usual, $Z(\g)$ denotes the centre of $U(\g)$. Standard properties of
the Harish-Chandra homomorphism $\varphi\colon\,Z(\g)\rightarrow
S(\h)$ ensure that the equality $Z(\g)\cap{\rm Ann}_{U(\g)}\,
L(\lambda')\,=\,Z(\g)\cap{\rm Ann}_{U(\g)}\, L(\lambda)$ holds if
and only if $\lambda'+\rho=w(\lambda+\rho)$ for some $w\in W$.

By Duflo's theorem [\cite{Du}, Theorem~1], for any primitive ideal
$\mathcal I$ of $U(\g)$ there exists an irreducible highest weight
module $L(\mu)$ such that ${\mathcal I}={\rm
Ann}_{U(\g)}\,L(\mu)$. Generically, the number of such
$\mu\in\h^*$ equals the order of $W$, but there are instances when
$\mu$ is uniquely determined by $\mathcal I$; this happens, for
example, when $\mathcal I$ has finite codimension in $U(\g)$. Note
that $Z(\g)\cap {\mathcal I}$ is a maximal ideal of $Z(\g)$ and
$z(\mu)=0$ for all $z\in\varphi(Z(\g)\cap \mathcal I)$. We say
that $\mathcal I$ has {\it rational infinitesimal character} if
$\langle\mu,\alpha^\vee\rangle\in{\mathbb Q}$ for all
$\alpha\in\Phi$ (our earlier remarks show that this is independent
of the choice of $\mu$). Our main goal in this note is to prove
the following\footnote{Recently I proved Theorem~\ref{Main} for
all primitive ideals $\mathcal I$ with $\mathcal{VA}({\mathcal
I})=\overline{\mathcal O}_\chi$ thereby confirming [\cite{P07},
Conjecture~3.2] in full generality. This result will appear
elsewhere.}:

\begin{theorem}\label{Main} If $\mathcal I$ is a primitive ideal of $U(\g)$ with
rational infinitesimal character and with $\mathcal{VA}(\mathcal
I)=\overline{\mathcal O}_\chi$, then ${\mathcal I}={\rm
Ann}_{U(\g)}(Q_\chi\otimes_{H_\chi}V)$ for some finite dimensional
irreducible $H_\chi$-module $V$.
\end{theorem}
\subsection{}\label{1.3}
To prove Theorem~\ref{Main} we work with a suitable localisation
$A=S^{-1}\Z$ of the ring of integers and consider natural
$A$-forms $\g_A$ and $L_A(\mu)$ of $\g$ and $L(\mu)$,
respectively. We have to choose $A$ very carefully; see Sect.~2.
After making our selection we pick $(p)\in{\rm Spec}\,A$ and
consider the highest weight module $L_p(\mu):=L_A(\mu)\otimes_A\k$
over the restricted Lie algebra $\g_\k=\g_A\otimes_A\k$, where
$\k$ is the algebraic closure of ${\mathbb F}_p$.

Let $G_\k$ be the simply connected algebraic $\k$-group with
$\Lie(G_\k)=\g_\k$. If $p\gg 0$, then $(\,\cdot\,,\,\cdot\,)$
induces a nondegenerate $G_\k$-invariant symmetric form on $\g_\k$
and there is a natural bijection between the nilpotent $G$-orbits
in $\g$ and the nilpotent $G_\k$-orbits in $\g_\k$. Let
$Z_p\,=\,\la x^p-x^{[p]}\,|\,\,x\in\g_\k\ra$, the $p$-centre of
the universal enveloping algebra $U(\g_\k)$. Given
$\eta\in\g_\k^*$ we denote by $\k_\eta=\k1_\eta$ the
$1$-dimensional $Z_p$-module such that
$(x^p-x^{[p]})(1_\eta)=\eta(x)^p1_\eta$ for all $x\in\g_\k$, and
we set $U_\eta(\g_\k):=\,U(\g_\k)\otimes_{Z_p}\k_\eta$ and
$L_p^\eta(\mu):=\,L_p(\mu)\otimes_{Z_p}\k_\eta$. The associative
algebra $U_\eta(\g_\k)$ is often referred to as the {\it reduced
enveloping algebra} of $\g_\k$ corresponding to $\eta$.

We show in (\ref{3.3}) that for all $p\gg 0$ the
$U_\eta(\g_\k)$-module $L_p^\eta(\mu)$ has dimension $\le
Dp^{d(e)}$, where $D=D(\mu)$ is independent of $p$ and
$d(e)=\frac{1}{2}\dim {\mathcal O}(\chi)$. Let ${\mathcal O}_\k$
be the nilpotent $G_\k$-orbit in $\g_\k$ corresponding to the
$G$-orbit containing $e$. We prove in (\ref{3.5}) that there is a
linear function $\bar{\chi}=(\bar{e},\,\cdot\,)$ on $\g_\k$ with
$\bar{e}\in{\mathcal O}_\k$ such that the
$U_{\bar{\chi}}(\g_\k)$-module $L_p^{\bar{\chi}}(\mu)$ is nonzero.
Here (and only here) we use our assumption that the infinitesimal
character of $\mathcal I$ is rational. Combining the
Kac--Weisfeiler conjecture proved in [\cite{P95}] with the above
results, we show in (\ref{3.6}) that the algebra
$U_{\bar{\chi}}(\g_\k)$ has a simple module of dimension
$kp^{d(e)}$ for some $k\le D=D(\mu)$. On the other hand, it is
established in [\cite{P02}] that
$U_{\bar{\chi}}(\g_\k)\cong\Mat_{p^{d(e)}}\big(H_{\bar{\chi}}^{[p]}\big)$,
where $H_{\bar{\chi}}^{[p]}$ is a `restricted' modular version of
the finite $W$-algebra $H_\chi\cong W^{\rm fin}(\g,e)$. We thus
deduce that for $p\gg 0$ the algebra $H_{\bar{\chi}}^{[p]}$ has a
simple module of dimension $k\le D=D(\mu)$.

The algebra $H_\chi$ relates to $H_{\bar{\chi}}^{[p]}$ in the same
way as $U(\g)$ relates to the restricted enveloping algebra
$U^{[p]}(\g_\k)$; in fact, when $\chi=0$ we have that
$H_\chi=U(\g)$ and $H_{\chi}^{[p]}=U^{[p]}(\g_\k)$. We use a PBW
basis of $H_\chi$ similar to the one described in (\ref{1.1}) to
find a suitable $A$-form $H_{\chi,\,A}$ in $H_\chi$ and show in
(\ref{4.5}) that $H_{\bar{\chi}}^{[p]}$ is a homomorphic image of
the $\k$-algebra $H_{\chi,\,\k}:=H_{\chi,\,A}\otimes_A\k$. Next we
introduce certain affine varieties ${\mathcal Y}_k$ of matrix
representations of $H_\chi$ and use reduction modulo $p$ together
with some results mentioned above to prove that ${\mathcal
Y}_l(\mathbb C)\ne\emptyset$ for some $l\le D$; see ({\ref{4.6}).
The definition of ${\mathcal Y}_l$ then implies that $H_\chi$
possesses an irreducible finite dimensional module $V$ with the
property that ${\mathcal I}\subseteq {\rm
Ann}_{U(\g)}(Q_\chi\otimes_{H_\chi} V)$. Since the primitive
ideals ${\rm Ann}_{U(\g)}(Q_\chi\otimes_{H_\chi} V)$ and
${\mathcal I}$ have the same associated variety, a well-known
result of Borho--Kraft [\cite{BK}] yields ${\mathcal I}={\rm
Ann}_{U(\g)}(Q_\chi\otimes_{H_\chi} V)$. Theorem~\ref{Main}
follows.
\begin{corollary}\label{AA} All finite
$W$-algebras $W^{\rm fin}(\g,e)$ possess finite dimensional
irreducible representations.
\end{corollary}
\noindent To deduce this corollary from Theorem~\ref{Main} we rely
in a crucial way on well-known results of Barbasch--Vogan
[\cite{BV1}, \cite{BV2}]. It would be interesting to find an
alternative proof.
\begin{corollary}\label{BB}
If $p\gg 0$, then for every linear function $\eta$ on $\g_\k$ the
reduced enveloping algebra $U_{\eta}(\g_\k)$ has a simple module of
dimension $lp^{\,(\dim G_\k\cdot\,\eta)/2}$ where $l<p$.
\end{corollary}
\noindent For $p$ large, Corollary~\ref{BB} provides a good upper
bound for the minimal dimension of irreducible
$U_\eta(\g_\k)$-modules and is the first result of its kind in the
literature. It can be regarded as a first step towards verifying
one of the main conjectures in the field which states, under mild
assumptions on $p$, that the minimal dimension of irreducible
$U_\eta(\g_\k)$-modules always equals $p^{\,(\dim
G_\k\cdot\,\eta)/2}$.

The general case of Corollary~\ref{BB} reduces quickly to the case
where $\eta=(\bar{e},\,\cdot\,)$ for some nilpotent element
$\bar{e}\in\g_\k$; see [\cite{P02}, 2.5] for more detail. In this
special case the the corollary follows from part~(2) of the proof
of Corollary~\ref{Walg} and the fact that $L_p^{\xi}(\mu)\ne 0$
for some $\xi\in G_\k\cdot \eta$; see Lemmas~\ref{lem1},
\ref{lem2} and Theorem~\ref{thm1}. We stress that this argument
also relies on the main results of [\cite{BV1}, \cite{BV2}]. It
would be very interesting to find a more direct argument.

\section{\bf Admissible rings for highest weight modules}
\subsection{}\label{2.1}
In what follows ${\mathbb Z}_+$ denotes the set of all nonnegative
integers. Given a commutative Noetherian ring $R$ we denote by
$\dim R$ the Krull dimension of $R$. Given a collection $\mathcal
C$ of regular functions on an algebraic variety $X$ we let
$V({\mathcal C})$ stand for the set of all common zeros of
$\mathcal C$ in $X$. If $\mathcal L$ is a Lie algebra over an
integral domain $A$, which is free as an $A$-module, then we write
$U_n({\mathcal L})$ for the $n$-th component of the canonical
filtration of the universal enveloping algebra $U(\mathcal L)$. By
the PBW theorem, the corresponding graded $A$-algebra $\gr\,
U(\mathcal L)$ is isomorphic to the symmetric algebra $S(\mathcal
L)$ of the free $A$-module $\mathcal L$. In particular, all
modules $U_n({\mathcal L})/U_{n-1}({\mathcal L})\cong S^n(\mathcal
L)$ are free over $A$.

Denote by ${\mathcal O}(x)$ the adjoint orbit $(\Ad G)\cdot x$ of
$x\in\g$ and put $d(x):=\frac{1}{2}\dim {\mathcal O}(x)$. We shall
often identify $\g$ with $\g^*$ via the Killing isomorphism
$\g\rightarrow \g^*$ which takes $x\in \g$ to the linear function
$(x,\,\cdot\,)$ on $\g$, where $(\,\cdot\,,\,\cdot\,)$ is a
multiple of the Killing form of $\g$ for which $(e,f)=1$.

Fix once and for all a primitive ideal ${\mathcal I}={\rm
Ann}_{U(\g)}\,L(\mu)$ such that $\mu(h_\alpha)\in\mathbb Q$ for all
$\alpha\in\Phi$ and $\mathcal{VA}({\mathcal I})=\overline{\mathcal
O}_\chi$. Write $\Phi^+=\{\gamma_1,\ldots,\gamma_N\}$, where $N$ is
the number of positive roots, put $\Phi^-:=-\Phi^+$, and let
$\g=\n^-\oplus\h\oplus\n^+$ be the corresponding triangular
decomposition of $\g$. Choose a Chevalley basis
$${\mathcal
B}=\{e_\gamma\,|\,\,\gamma\in\Phi\}\cup\{h_\alpha\,|\,\,\alpha\in\Pi\}$$
in $\g$ and set ${\mathcal
B}^{\pm}:=\{e_\alpha\,|\,\,\alpha\in\pm\Phi^+\}$. Given a (unitary)
subring $A$ of $\mathbb C$ we denote by $\g_A$ and $\n_A^\pm$ the
$A$-spans of $\mathcal B$ and ${\mathcal B}^\pm$, respectively.
Denote by $U(\g_A)$ and $U(\n_A^\pm)$ the $A$-lattices in $U(\g)$
and $U(\n^\pm)$ generated as $A$-algebras by $\mathcal B$ and
${\mathcal B}^-$, respectively.

We call a subring $A$ {\it admissible} if $A=S^{-1}\mathbb Z$ where
$S$ is a finitely generated multiplicative subset of $\mathbb Z$
containing denominators of all $\mu(h_\alpha)$ with $\alpha\in\Pi$.
Note that $A$ is a Noetherian principal ideal domain, hence a
Dedekind ring. Denote by $\pi(A)$ the set of all primes $p\in
\mathbb Z$ with $pA\ne A$. As $S$ is finitely generated, $\pi(A)$
contains almost all primes in $\mathbb Z$.

The Killing form $\kappa$ of $\g$ is $\Z$-valued on $\g_{\Z}$ and
$\kappa(e,f)\in \Z\setminus\{0\}$ by the $\sl_2$-theory. {\it We
always choose $A$ such that $\kappa(e,f)$ is invertible in $A$}.
This will ensure that the form $(\,\cdot\,,\,\cdot\,)$ is
$A$-valued on $\g_A$. {\it We always choose $A$ such that the
determinant of the Gram matrix of $(\,\cdot\,,\,\cdot\,)$ relative
to the Chevalley basis $\mathcal B$ and all bad primes of the root
system $\Phi$ are invertible in $A$}.
\subsection{}\label{2.2}
Let $M(\mu)$ be the Verma module of highest weight $\mu$, and let
$\tilde{v}_\mu$ be a highest weight vector of weight $\mu $ in
$M(\mu)$. Given a subring $A$ of $\mathbb Q$ we denote by $M_A(\mu)$
the $A$-span of all $e_{\gamma_1}^{a_1}\cdots
e_{\gamma_N}^{a_N}\tilde{v}_\mu$ with $a_i\in{\mathbb Z}_+$. These
vectors form a free basis of the $A$-module $M_A(\mu)$, so that
$M_A(\mu)=U(\n_A^-)v_0\cong U(\n_A^-)$ as $A$-modules. Because $A$
is admissible, $M_A(\mu)$ is a $\g_A$-stable $A$-lattice in
$M(\mu)$.

Given an $\h$-module $V$ we denote by ${\mathcal X}(V)$ the set of
all weights of $V$ relative to $\h$ and write $V_\nu$ for the weight
space of $V$ corresponding to $\nu\in {\mathcal X}(V)$. If $E$ is an
$A$-submodule of $V$, then we set $E_\nu:=E\cap V_\nu$. By the PBW
theorem,
$$M_A(\mu)\ =\bigoplus_{\nu\in\, {\mathcal X}(M(\mu))}\,M_A(\mu)_\nu \ \ \mbox{ and
}\ \ S(\n_A^-)\ =\bigoplus_{\nu\in\, {\mathcal
X}(S(\n^-))}\,S(\n_A^-)_\nu.$$  All weight components of $M_A(\mu)$
and $S(\n_A^-)$ are free $A$-modules of finite rank. Denote by
$M^{\rm max}(\mu)$ the unique maximal submodule of $M(\mu)$, so that
$L(\mu)=M(\mu)/M^{\rm max}(\mu)$, and let $v_\mu$ denote the image
of $\tilde{v}_\mu$ under the canonical homomorphism
$M(\mu)\twoheadrightarrow L(\mu)$. Put $M_A^{\rm max}(\mu):=M^{\rm
max}(\mu)\cap M_A(\mu)$ and define
$$L_A(\mu)\,:=\,M_A(\mu)/M_A^{\rm max}(\mu).$$

For $n\in\Z_+$ we define $L_n(\mu):=U_n(\g)v_\mu=U_n(\n^-)v_\mu$ and
$L_{A,n}(\mu):=U_n(\g_A)v_\mu=U_n(\n_A^-)v_\mu$, and set
$$\gr\,L(\mu)=\,\textstyle{\bigoplus_{n\ge
0}}\,\,L_n(\mu)/L_{n-1}(\mu)\quad\mbox{and}\quad
\gr\,L_A(\mu)=\,\textstyle{\bigoplus_{n\ge
0}}\,\,L_{A,n}(\mu)/L_{A,n-1}(\mu)$$ (here $L_{-1}(\mu)=0$ and
$L_{A,-1}(\mu)=0$ by convention). Note that $\gr\,L(\mu)$ and
$\gr\,L_A(\mu)$ are generated by $v_\mu=\gr_0\,v_\mu$ as modules
over $S(\g)=\gr\,U(\g)$ and $S(\g_A)=\gr\, U(\g_A)$, respectively.
We now define
$${\mathcal J}:={\rm Ann}_{S(\g)}\,\gr\,L(\mu)\,=\,{\rm
Ann}_{S(\g)}\,v_\mu\quad \mbox{and}\quad {\mathcal J}_A:={\rm
Ann}_{S(\g_A)}\,\gr\,L(\mu)\,=\,{\rm Ann}_{S(\g_A)}\,v_\mu.$$ These
are graded ideals of $S(\g)$ and $S(\g_A)$, respectively. We put
${\mathcal R}:=S(\g)/{\mathcal J}$ and ${\mathcal
R}_A:=S(\g_A)/{\mathcal J}_A$.  It is worth remarking that
$L_{\mathbb Q}(\mu)$ and ${\mathcal R}_{\mathbb Q}$ are ${\mathbb
Q}$-forms in $L(\mu)$ and $\mathcal R$, respectively.

The zero locus of the ideal ${\mathcal J}$ in $\g^*={\rm
Specm}\,S(\g)$ is called the {\it associated variety} of the
$\g$-module $L(\mu)$ and denoted by ${\mathcal V}_\g L(\mu)$. It
follows from a theorem of Gabber that all irreducible components of
the variety ${\mathcal V}_\g L(\mu)$ have dimension $d(e)$; see
[\cite{J84}] for more detail. In particular, $\dim\,{\mathcal
R}=d(e)$.

Suppose now that $A$ is admissible. Since $M^{\rm max}(\mu)$ is the
kernel of the Shapovalov form of $M(\mu)$, it is straightforward to
see that $M_A^{\rm max}(\mu)$ is a direct sum of its weight
components $M_A^{\rm max}(\mu)_\nu$ and each weight component
$L_A(\mu)_\nu=\,M_A(\mu)_\nu/M_A^{\rm max}(\mu)_\nu\,$ of $L_A(\mu)$
is finitely generated and torsion free as an $A$-module. Since $A$
is a principal ideal domain, it follows that all $L_A(\mu)_\nu$ are
free $A$-modules of finite rank; see [\cite{Bo}, Ch.~VII,
Sect.~4.10].

By construction, ${\mathcal R}_A\,\cong\,\gr\,L_A(\mu)$ as graded
$A$-modules. An element $x\in {\mathcal R}_A$ is said to be a {\it
torsion element} if there is a nonzero $a\in A$ such that $a x=0$.
It is easy to see that the set ${\mathcal R}_A^{\rm tor}$ of all
torsion elements of ${\mathcal R}_A$ is an ideal of ${\mathcal
R}_A$. Since $A$ is admissible, ${\mathcal R}_A$ is a finitely
generated $\mathbb Z$-algebra, hence a commutative Noetherian
ring. Therefore, the ideal ${\mathcal R}_A^{\rm tor}$ is a
finitely generated ${\mathcal R}_A$-module. So there exists an
integer $b$ such that $bx=0$ for {\it all} $x\in{\mathcal
R}_A^{\rm tor}$. It follows that ${\mathcal R}_A\otimes_A
A[b^{-1}]$ is torsion free over $A[b^{-1}]$. On the other hand,
${\mathcal J}_{A[b^{-1}]}={\mathcal J}\otimes_A A[b^{-1}]$ and
$${\mathcal R}_{A[b^{-1}]}=S(\g_{A[b^{-1}]})/{\mathcal J}_{A[b^{-1}]}\cong
\big(S(\g_A)\otimes_A A[b^{-1}]\big)/\big({\mathcal
J}_A\otimes_AA[b^{-1}]\big)\cong{\mathcal R}_A\otimes_A A[b^{-1}]$$
by standard properties of localisation; see [\cite{Bo}, Ch.~II,
2.4]. Thus, the $A[b^{-1}]$-module ${\mathcal
R}_{A[b^{-1}]}\cong\gr\,L_{A[b^{-1}]}(\mu)$ is torsion free.

{\it From now on we always assume that ${\mathcal R}_A$ is a
torsion free $A$-module}. This assumption is justified by the
discussion above.
\subsection{}\label{2.3}
Since $A$ is a principal ideal domain and ${\mathcal R}_A\cong
\gr\,L_A(\mu)$ is torsion free over $A$, all graded components
${\mathcal R}_{A,n}\cong L_{A,n}(\mu)/L_{A,n-1}(\mu)$ of ${\mathcal
R}_A$ are free $A$-modules of finite rank. Moreover, ${\rm
rk}_A\,{\mathcal R}_{A,n}=\,\dim_{\mathbb C}{\mathcal
R}_n=\dim_{\mathbb C}L_n(\mu)/L_{n-1}(\mu)$ for all $n$. It follows
that ${\mathcal R}_{\mathbb Q}\,\cong\,{\mathcal
R}_{A}\otimes_A\mathbb Q$ is a $\mathbb Q$-form of $\mathcal R$.
Since ${\mathcal R}\,\cong {\mathcal R}_{\mathbb Q}\otimes_{\mathbb
Q}{\mathbb C}$ as algebras over $\mathbb C$, it is immediate from
[\cite{Eis}, Corollary~13.7] that $\dim {\mathcal R}_{\mathbb
Q}=\dim{\mathcal R}=\,d(e)$.

According to Noether Normalisation, there exist homogeneous,
algebraically independent $y_1,\ldots, y_d\in{\mathcal R}_{\mathbb
Q}$, where $d=d(e)$, such that ${\mathcal R}$ is a finitely
generated module over its graded polynomial subalgebra ${\mathbb
Q}[y_1,\ldots, y_d]$; see [\cite{Eis}, Theorem~13.3] (the elements
$y_1,\ldots, y_d$ remain algebraically independent in ${\mathcal
R}\,\cong {\mathcal R}_{\mathbb Q}\otimes_{\mathbb Q}{\mathbb C}$).
Let $v_1,\ldots, v_D$ be a generating set of the ${\mathbb
Q}[y_1,\ldots, y_d]$-module ${\mathcal R}_{\mathbb Q}$ and let
$x_1,\ldots, x_{m'}$ be a generating set of the $A$-algebra
${\mathcal R}$. Then
\begin{eqnarray*}
v_i\cdot v_j&=&\textstyle{\sum}_{k=1}^D \,p_{i,j}^k(y_1,\ldots,
y_d)v_k\ \ \qquad(1\le
i,j\le D)\\
x_i&=&\textstyle{\sum}_{j=1}^{D}\,q_{i,j}(y_1,\ldots, y_d) v_j\ \
\,\qquad(1\le i\le m')
\end{eqnarray*} for some polynomials $p_{i,j}^k,\, q_{i,j}\in
{\mathbb Q}[X_1,\ldots, X_d].$ Recall that ${\mathcal R}_A$ contains
a basis of ${\mathcal R}_{\mathbb Q}$ over $\mathbb Q$. The
coordinate vectors of the $x_i$'s, $y_i$'s and $v_i$'s relative to
this basis and the coefficients of the polynomials $q_{i,j}$ and
$p_{i,j}^k$ involve only finitely many scalars in $\mathbb Q$.
Enlarging our multiplicative set $S$ if need be we may assume that
all $y_i$ and $v_i$  are in ${\mathcal R}_A$ and all $p_{i,j}^k$ and
$q_{i,j}$ are in $A[X_1,\ldots,X_d]$. {\it Thus, we may assume that
$${\mathcal R}_A\,=\,A[y_1,\ldots,
y_d]v_1+\cdots+A[y_1,\ldots,y_d]v_D$$ is a finitely generated module
over the polynomial algebra $A[y_1,\ldots, y_d]$}.
\subsection{}\label{2.4}
We let the Galois group $\Gamma={\rm Gal}({\mathbb
C}/{\mathbb Q})$ act on $L(\mu)$ by semilinear automorphisms fixing
all elements of the $\mathbb Q$-form $L_{\mathbb Q}(\mu)$. Note that
$\Gamma$ acts on $U(\g)$ as semilinear algebra automorphisms fixing
all elements of the $\mathbb Q$-form $U(\g_{\mathbb Q})$ of $U(\g)$.
These two Galois actions are compatible, that is $\sigma(u.
v)=u^\sigma.\,\sigma(v)$ for all $u\in U(\g)$, $v\in L(\mu)$ and
$\sigma\in\Gamma$. Therefore, ${\mathcal I}^\sigma={\mathcal I}$ for
all $\sigma\in\Gamma$. Since $\Gamma$ preserves the canonical
filtration of $U(\g)$, each subspace ${\mathcal I}\cap U_n(\g)$ is
$\Gamma$-invariant. It follows that the graded ideal
$$\gr\,{\mathcal I}\,=\,\textstyle{\bigoplus}_{n\ge
0}\,\big({\mathcal I}\cap U_n(\g)\big)/\big({\mathcal I}\cap
U_{n-1}(\g)\big)\subset S(\g)$$ is $\Gamma$-stable. As a
consequence, for every $n\in\mathbb Z_+$ the $\mathbb Q$-subspace
$S^n(\g_{\mathbb Q})\cap\gr\,{\mathcal I}$ is a $\mathbb Q$-form of
the graded component $\gr_n\,{\mathcal I}\subset S^n(\g)$. Since
$S(\g)$ is Noetherian, the ideal $\gr\,\mathcal I$ is generated over
$S(\g)$ by its $\mathbb Q$-subspace $\gr\,{\mathcal I}_{{\mathbb
Q},n'}:=\big(\bigoplus_{k\le n'}S^k(\g_{\mathbb
Q})\big)\cap\gr\,{\mathcal I}$ for some $n'=n'(\mu)\in\Z_+$. Using
standard filtered-graded techniques it is now easy to derive that
$\mathcal I$ is generated over $U(\g)$ by its $\mathbb Q$-subspace
${\mathcal I}_{{\mathbb Q},n'}:=\,U_{n'}(\g_{\mathbb Q})\cap\mathcal
I$.

Since $\mathcal I$ is a two-sided ideal of $U(\g)$, all subspaces
${\mathcal I}\cap U_n(\g)$ and $\gr_n\,\mathcal I$ are invariant
under the adjoint action of $G$ on $U(\g)$. It follows that the
$\mathbb Q$-subspaces $\gr\,{\mathcal I}_{{\mathbb Q},n'}$ and
${\mathcal J}_{{\mathbb Q},n'}$ are stable under the adjoint action
of the Kostant $\Z$-form $U_\Z$ of $U(\g)$. Since $\h_{\mathbb
Q}:=\h\cap\g_{\mathbb Q}$ is a split Cartan subalgebra of
$\g_{\mathbb Q}$, the adjoint $\g_{\mathbb Q}$-modules
$\gr\,{\mathcal I}_{{\mathbb Q},n'}$ and ${\mathcal I}_{{\mathbb
Q},n'}$ decompose into a direct sum of absolutely irreducible
$\g_{\mathbb Q}$-modules with integral dominant highest weights.
Consequently, these $\g_{\mathbb Q}$-modules possess $\Z$-forms
invariant under the adjoint action of $U_{\mathbb Z}$; we call them
$\gr\,{\mathcal I}_{\Z, n'}$ and ${\mathcal I}_{\Z,n'}$.

It is immediate from Weyl's theorem that the $\g_\mathbb
Q$-modules $\gr\,{\mathcal I}_{{\mathbb Q},n'}$ and ${\mathcal
I}_{{\mathbb Q},n'}$ are isomorphic. Let $\{\psi_i\,|\,\,i\in I\}$
be a homogeneous basis of the free $\Z$-module $\gr\,{\mathcal
I}_{\Z, n'}$ and let $\{u_i\,|\,\,i\in I\}$ be any basis of the
free $\Z$-module ${\mathcal I}_{\Z,n'}$.  Expressing the $u_i$ and
$\psi_i$ via the PBW bases of $U(\g_{\mathbb Q})$ and
$S(\g_{\mathbb Q})$ associated with our Chevalley basis $\mathcal
B$ brings up only finitely many scalars in $\mathbb Q$. Thus, no
generality will be lost by assuming that all these scalars are in
$A$. In other words, {\it it can be assumed that all $\psi_i$ are
in $S(\g_A)$ and all $u_i$ are in $U(\g_A)$}.
\subsection{}\label{2.5} Let $K$ be an algebraically closed field
whose characteristic is a good prime for the root system $\Phi$. Let
$\g_K=\g_{\Z}\otimes_{\Z}K$ and let $G_K$ be the simple, simply
connected algebraic $K$-group with hyperalgebra
$U_K:=U_{\Z}\otimes_{\Z}K$. Let ${\mathcal N}(\g)$ and ${\mathcal
N}(\g_K)$ denote the nilpotent cones of $\g$ and $\g_K$,
respectively. It is well-known that there is a natural bijection
between the $G$-orbits in ${\mathcal N}(\g)$ and the $G_K$-orbits in
${\mathcal N}(\g_K)$. More precisely, combining [\cite{SpSt}, III,
4.29] with [\cite{P03}, Theorems~2.6 and 2.7] it is not hard to
observe that there are nilpotent elements $e_1,\ldots,
e_t\in\g_{\Z}$ such that
\begin{itemize}
\item[(i)\,] $\{e_1,\ldots, e_t\}$ is a set of representatives for
$\,{\mathcal N}(\g)/G$;

\smallskip

\item[(ii)\,] $\{e_1\otimes 1,\ldots, e_t\otimes 1\}$ is a set of
representatives for $\,{\mathcal N}(\g_K)/G_K$;

\smallskip

\item[(iii)\,] $\,\dim_{\mathbb C}\,(\Ad G)e_i\,=\,\dim_K\,(\Ad
G_K)(e_i\otimes 1)\, $ for all $\,i\le t$.
\end{itemize}
No generality will be lost by assuming that $e=e_k$ for some $k\le
t$ and ${\mathcal O}(e_i)\subset \overline{{\mathcal O}(e)}$ for all
$i\le k$.

For $1\le i\le t$ we set $\chi_i:=(e_i,\,\cdot\,).$ Since
$\mathcal{VA}({\mathcal I})$ is the zero locus of $\gr\, {\mathcal
I}$ and the ideal $\gr\,\mathcal I$ is generated by the set
$\{\psi_i\,|\,\,i\in I\}$, we have that $\overline{{\mathcal
O}(\chi)}=V(\gr\,{\mathcal I})=\bigcap_{i\in I}\,V(\psi_i)$. It
follows that the $\psi_i$ vanish on all $\chi_j$ with $j\le k$.
Since all $\psi_i$ are in $S(\g_A)$, all $e_j$ are in $\g_{\Z}$, and
the form $(\,\cdot\,,\,\cdot\,)$ is $A$-valued, we also have that
$\psi_i(\chi_j)\in A$. We now impose our next assumption on $A$:
{\it all nonzero $\psi_i(\chi_j)$ are invertible in $A$}.
\section{\bf Highest weight modules and non-restricted representations}
\subsection{}\label{3.1} In this section we assume that our admissible ring $A=S^{-1}\Z$
for the $\g$-module $L(\mu)$ satisfies all the requirements
highlighted in (\ref{2.1})\,--\,(\ref{2.5}). We fix $p\in\pi(A)$
and set $L_p(\mu):=L_A(\mu)\otimes_{A}\k$ where $\k$ is the
algebraic closure of ${\mathbb F}_p$. Let
$\bar{v}_\mu:=v_\mu\otimes 1$, the image of $v_\mu$ in $L_p(\mu)$.
Clearly, $L_p(\mu)$ is a module over the Lie algebra
$\g_{\k}=\g_{A}\otimes_A\k$. Given $x\in\g_A$ we denote by
$\bar{x}$ the image of $x$ in $\g_\k$, that is $\bar{x}=x\otimes
1$. It is immediate from our assumptions on $A$ made in
(\ref{2.1}) that $\g_{\k}$ is a simple Lie algebra with $\k$-basis
$\bar{\mathcal B}:=\{\bar{e}_\gamma\,|\,\,\gamma\in\Phi\}
\cup\{\bar{h}_\alpha\,|\,\alpha\in\Pi\}$.  Since the Lie algebra
$\g_\k$ is simple and restricted, it carries a unique $p$-mapping
$x\mapsto x^{[p]}$ which has the property that
$\bar{e}_\gamma^{[p]}=0$ for all $\gamma\in\Phi$ and
$\bar{h}_\alpha^{[p]}=\bar{h}_\alpha$ for all $\alpha\in\Pi$.

The $A$-valued bilinear form $(\,\cdot\,,\,\cdot\,)$ on $\g_A$
gives rise to a $\k$-bilinear form on $\g_k$ which, for ease of
notation, shall be denoted by the same symbol. Our assumptions in
(\ref{2.1}) ensure that $(\,\cdot\,,\,\cdot\,)$ is nondegenerate
on $\g_\k$. Let $G_\k$ denote the simple, simply connected
algebraic $\k$-group with hyperalgebra
$U_\k:=U_{\Z}\otimes_{\Z}\k$. By construction,
$\g_\k\,=\,\Lie(G_\k)$; so the group $G_\k$ as on $\g_\k$ as Lie
algebra automorphisms. It is well-known (and straightforward to
see) that the form $(\,\cdot\,,\,\cdot\,)$ is $G_\k$-invariant. We
shall often identify the $G_\k$ modules $\g_\k$ and $\g_\k^*$ by
using this nondegenerate bilinear form.

Now $\g_\k=\n^{-}_\k\oplus\h_\k\oplus\n^{+}_\k$, where
$\n^{\pm}_\k:=\n^{\pm}\otimes_A\k$ and $\h_\k:=\h_A\otimes_A\k$.
Furthermore, $\bar{v}_\mu$ is a weight vector for the torus $\h_\k$
of the restricted Lie algebra $\g_\k$, the subalgebra $\n^+_\k$
annihilates $\bar{v}_\mu$, and $L_p(\mu)= U(\n^-_\k)\cdot
\bar{v}_\mu$. We denote by $\bar{\mu}$ the $\h_\k$-weight of
$\bar{v}_\mu$. Since $\mu(h_\alpha)\in A$ for all $\alpha\in \Pi$
and $A\otimes_{A}\k={\mathbb F}_p$, we have that
$\bar{\mu}(\bar{h}_\alpha)\in{\mathbb F}_p$ for all $\alpha\in \Pi$.

Let $Z_p$ denote the $p$-centre of the universal enveloping
algebra $U(\g_\k)$, the central subalgebra of $U(\g_\k)$ generated
by all $x^p-x^{[p]}$ with $x\in\g_\k$. The PBW theorem implies
that $Z_p$ is a polynomial algebra in $\bar{e}_\gamma^p$ with
$\gamma\in\Phi$ and $\bar{h}_\alpha^p-\bar{h}_\alpha$ with
$\alpha\in\Pi$, and $U(\g_\k)$ is a free $Z_p$-module of rank
$p^{\dim\g}$. Given a linear function $\xi$ on $\g_\k^*$ we denote
by $I_\xi$ the two-sided ideal of $U(\g_\k)$ generated by the
central elements $x^p-x^{[p]}-\xi(x)^p1$, where $x\in\g_\k$, and
set $U_\xi(\g_\k):=U(\g_k)/I_\xi$. The associative algebra
$U_\xi(\g_\k)$ is called the {\it reduced enveloping algebra} of
$\g_\k$ associated with $\xi$. The above remark shows that
$\dim_\k U_\xi(\g_\k)=p^{\dim\g}$.

Every irreducible $\g_\k$-module is a module over $U_\xi(\g_\k)$ for
precisely one $\xi=\xi_V\in\g_\k^*$. The linear function $\xi_V$ is
called the {\it $p$-character} of $V$; see [\cite{P95}] for more
detail.

\subsection{}\label{3.2}
For every maximal ideal $J$ of $Z_p$ there exists a unique linear
function $\eta=\eta_J$ on $\g_\k$ such that $J$ is generated by
all $x^p-x^{[p]}-\eta(x)^p1$ with $x\in\g_\k$, so that $J=Z_p\cap
I_\eta$. As the Frobenius map of $\k$ is bijective, this allows us
to identify ${\rm Specm}\,Z_p$ with $\g_\k^*$.

We now define $I_p(\mu):=\{z\in Z_p\,|\,\,z\cdot \bar{v}_\mu=0\}$,
an ideal of the $p$-centre $Z_p$, and denote by $V_p(\mu)$ the zero
locus of $I_p(\mu)$ in $\g_\k^*$. Our remarks in (\ref{3.1}) show
that $\bar{e}_\gamma^p\in I_p(\mu)$ for all $\gamma\in\Phi^+$ and
${\bar{h}^p_\alpha}-\bar{h}_\alpha\in I_p(\mu)$ for all
$\alpha\in\Pi$ (it is important here that
$\bar{\mu}(\bar{h}_\alpha)\in\mathbb F_p$ for all $\alpha\in\Pi$).
It follows that
\begin{eqnarray}\label{support} V_p(\mu)\,\subseteq\,
\{\eta\in\g_\k^*\,|\,\,\,\eta(\h_\k)=\eta(\n^+_\k)=0\}.
\end{eqnarray} Since the $U(\g_\k)$-module $L_p(\mu)$ is generated by
$\bar{v}_\mu$, we have that $I_p(\mu)\,=\,{\rm
Ann}_{Z_p}\,L_p(\mu)$. Given $\eta\in\g_\k^*$ we set
$L_p^\eta(\mu):=L_p(\mu)/I_\eta\cdot L_p(\mu).$ By construction,
$L_p^\eta(\mu)$ is a $\g_\k$-module with $p$-character $\eta$.
\begin{lemma}\label{lem1}
If $\eta\in V_p(\mu)$, then $L_p^\eta(\mu)$ is a nonzero
$U_\eta(\g_\k)$-module.
\end{lemma}\label{lem2}
\begin{proof}
Suppose $L_p^\eta(\mu)=0$ for some $\eta\in V_p(\mu)$, and set
$J_\eta:=Z_p\cap I_\eta$. Then $L_p(\mu)=I_\eta\cdot
L_p(\mu)=J_\eta\cdot L_p(\mu)$. It follows from our discussion in
(\ref{3.1}) that $L_p(\mu)$ is a finitely generated $Z_p$-module. By
a version of Nakayama's lemma, this implies that there is $z\in
J_\eta$ such that $1-z\in {\rm Ann}_{Z_p}\,L_p(\mu)$; see
[\cite{Eis}, Corollary~4.7] for example. On the other hand, ${\rm
Ann}_{Z_p}\,L_p(\mu)=I_p(\mu)\subseteq \sqrt{I_p(\mu)}\subseteq
J_\eta$ because $\eta\in V_p(\mu)$. But then both $z$ and $1-z$ are
in $J_\eta$, a contradiction.
\end{proof}
\begin{rem}\label{rem1} It follows from (\ref{support}) that every linear function $\xi\in
V_p(\mu)$ is {\it nilpotent}. More precisely, $\xi=(x,\cdot\,)$ for
some $x\in\n^+_\k$.
\end{rem}
\subsection{}\label{3.3}
Set $L_{p,n}(\mu):=U_n(\g_\k)\bar{v}_\mu$ and
$\gr\,L_p(\mu):=\,\textstyle{\bigoplus_{n\ge
0}}\,\,L_{p,n}(\mu)/L_{p, n-1}(\mu),$ where $n\in\Z_+$. Note that
$\gr\,L_p(\mu)$ is a cyclic $S(\g_\k)$-module generated by
$\bar{v}_\mu=\gr_0\, \bar{v}_\mu$. Also,
$$L_{p,n}(\mu)=U_n(\g_k)\bar{v}_\mu=\big(U_n(\g_A)v_\mu\big)\otimes_A\k
=L_{A,n}(\mu)\otimes_A\k.$$ We put ${\mathcal J}_p:={\rm
Ann}_{S(\g_\k)}\,\gr\,L_p(\mu)={\rm Ann}_{S(\g_\k)}\,\bar{v}_\mu$
and ${\mathcal R}_p:=S(\g_\k)/{\mathcal J}_p$, and denote by
${\mathcal V}_\g L_p(\mu)$ the zero locus of ${\mathcal J}_p$ in
${\rm Specm}\,S(\g_\k)=\g_\k^*$. Since $\bar{v}_\mu$ is a highest
weight vector for $\h_\k\oplus\n_\k^+$, all linear functions from
${\mathcal V}_\g L_p(\mu)$ vanish on $\h_\k\oplus\n_\k^+$.

By our assumptions in (\ref{2.2}) and (\ref{2.3}), all $L_{A,n}$ and
all graded components ${\mathcal R}_{A,n}\cong
L_{A,n}(\mu)/L_{A,n-1}(\mu)$ of ${\mathcal R}_A$ are free
$A$-modules of finite rank. From this it is immediate that
${\mathcal R}_p\cong {\mathcal R}_A\otimes_A\k$ as graded
$\k$-algebras. Comparing the Hilbert polynomials of ${\mathcal
R}={\mathcal R}_A\otimes_A\mathbb C$ and ${\mathcal R}_p$ we now
deduce that $\dim {\mathcal R}_p\,=\,\dim{\mathcal R}\,=\,d(e);$ see
[\cite{Eis}, Corollary~13.7]. As a consequence,
\begin{eqnarray}\label{Krull}\dim_\k {\mathcal V}_\g L_p(\mu)\,=\,
\dim {\mathcal R}_p=\,d(e).
\end{eqnarray}
Recall that ${\mathcal R}_A=A[y_1,\ldots,y_d]v_1+\cdots
+A[y_1,\ldots, y_d]v_D,$ for some $v_1,\ldots, v_D\in {\mathcal
R}_A$, where $d=d(e)$ and $A[y_1,\ldots, y_d]$ is a polynomial
algebra over $A$ contained in ${\mathcal R}_A$; see our discussion
in (\ref{2.3}). We stress that $D=D(\mu)$ depends on $\mu$, but not
on $p$.
\begin{lemma}
For every $\eta\in\g_\k^*$ we have that $\,\dim_\k L_p^\eta(\mu)\le
Dp^{d(e)}$.
\end{lemma}
\begin{proof}
Suppose that $y_i$ has degree $a_i$, where $1\le i\le d$, and
$v_k$ has degree $l_k$, where $1\le k\le D$, and let
$\Phi_A\colon\,S(\g_A)\twoheadrightarrow {\mathcal R}_A$ denote
the canonical homomorphism. For $1\le i\le d$ (resp. $1\le k\le
D$) choose $u_i\in U_{a_i}(\g_A)$ (resp. $w_k\in U_{l_k}(\g_A)$)
such that $\Phi_A(\gr_{a_i}\,u_i)=y_i$ (resp.
$\Phi_A(\gr_{l_k}\,w_k)=v_k$). Let $\bar{u}_i$ (resp. $\bar{w}_k$)
denote the image of $u_i$ (resp. $w_k$) in
$U(\g_\k)=U(\g_A)\otimes_A\k$. For each $n\in\Z_+$ the set
$$\{w_ku_1^{i_1}\cdots u_d^{i_d}\cdot v_{\mu}\,|\,\,
l_k+\textstyle{\sum}_{j=1}^d\, i_ja_j\le n;\ 1\le k\le D\}$$ spans
the $A$-module $L_{A,n}(\mu)$. In view of our earlier remarks this
implies that the set
$$\{\bar{w}_k\bar{u}_1^{i_1}\cdots \bar{u}_d^{i_d}\cdot \bar{v}_{\mu}\,|\,\,
l_k+\textstyle{\sum}_{j=1}^d\, i_ja_j\le n;\ 1\le k\le D\}$$ spans
the $\k$-space $L_{p,n}(\mu)$. Since
$\gr_{pa_i}(\bar{u}_i^p)=(\gr_{a_i}\bar{u}_i)^p$ is a $p$-th power
in $S(\g_\k)$, for every $i\le d$ there exists a $z_i\in Z_p\cap
U_{a_i}(\g_\k)$ such that $\bar{u}_i^p-z_i\in U_{pa_i-1}(\g_\k)$.

We denote by $L_p'(\mu)$ the $Z_p$-submodule of $L_p(\mu)$ generated
by all $\bar{w}_k\bar{u}_1^{i_1}\cdots \bar{u}_d^{i_d}\cdot
\bar{v}_{\mu}$ with $0\le i_j\le p-1$ and  $1\le k\le D.$ Using the
preceding remarks and induction on $n$ we now obtain that
$L_{p,n}(\mu)\subset L_p'(\mu)$ for all $n\in\Z_+$. But then
$L_p(\mu)=L_p'(\mu)$, implying that the image of
$\{\bar{w}_k\bar{u}_1^{i_1}\cdots \bar{u}_d^{i_d}\cdot
\bar{v}_{\mu}\,|\,\, 0\le i_j\le p-1;\, 1\le k\le D\}$ in
$L_p^\eta(\mu)=L_p(\mu)/I_\eta\cdot L_p(\mu)$ spans $L_p^\eta(\mu)$
for every $\eta\in\g_\k^*$.  Since $d=d(e)$, this yields $\dim_\k
L_p^\eta(\mu)\le Dp^{d(e)}$ completing the proof.
\end{proof}
\subsection{}\label{3.4} {\it From now on we shall always assume that $D!$
is invertible in $A$}. Then $p>D$ for every $p\in\pi(A)$. We shall
identify $\g_\k$ with $\g_\k^*$ by using the $G_\k$-equivariant
map $x\mapsto (x,\cdot\,)$. Under this identification, we have
that $V_p(\mu)\subseteq \n^+_\k$; see Remark~\ref{rem1}. The
$p$-centre $Z_p(\n_\k^-)=Z_p\cap U(\n_\k^-)$ of the enveloping
algebra $U(\n_\k^-)$ is isomorphic to the polynomial algebra
$\k[\bar{e}_{-\gamma_1}^p\,\ldots,\bar{e}_{-\gamma_N}^p]$, hence
can be canonically identified with the subalgebra $S(\n_\k^-)^p$
of all $p$-th powers in $S(\n_\k^-)$. Thus, we can regard
$I_p(\mu)\cap Z_p(\n_\k^-)$ as an ideal of the graded polynomial
algebra
$S(\n_\k^-)^p=\k[\bar{e}_{-\gamma_1}^p\,\ldots,\bar{e}_{-\gamma_N}^p]$.
It follows from our discussion in (\ref{3.2}) and the above
identifications that
\begin{eqnarray}\label{VIP}
V(I_p(\mu)\cap Z_p(\n_\k^-))\cap\n_\k^+=V_p(\mu).
\end{eqnarray}
Let $\gr\big(I_p(\mu)\cap Z_p(\n_\k^-)\big)$ be the homogeneous
ideal of $S(\n_\k^-)^p$ spanned by the highest components of all
elements in $I_p(\mu)\cap Z_p(\n_\k^-)$. From (\ref{VIP}) it
follows that the zero locus of $\gr\big(I_p(\mu)\cap
Z_p(\n_\k^-)\big)$ in $\n_\k^+$ coincides with ${\mathbb
K}(V_p(\mu))$, the associated cone to $V_p(\mu)$; see [\cite{Kr},
II, 4.2]. Since $I_p(\mu)\cap Z_p(\n_\k^-)$ is contained in ${\rm
Ann}_{Z_p}\,\bar{v}_\mu$, all elements of $\gr\big(I_p(\mu)\cap
Z_p(\n_\k^-)\big)$ annihilate $\gr_0\,
\bar{v}_\mu\in\gr\,L_p(\mu)$. Then $\gr\big(I_p(\mu)\cap
Z_p(\n_\k^-)\big)\subset {\mathcal J}_p\cap S(\n_\k^-)$, which
yields
\begin{eqnarray}\label{VA-VP}
{\mathcal V}_\g L_p(\mu)\,=V({\mathcal J}_p\cap
S(\n_\k^-))\cap\n_\k^+\,\subseteq\, {\mathbb K}(V_p(\mu)).
\end{eqnarray}
\subsection{}\label{3.5} Given $x\in\g_\k$ denote by ${\mathcal
O}_\k(x)$ the adjoint $G_\k$-orbit containing $x$. For $1\le i\le
t$ set $\bar{e}_i=e_i\otimes 1$. By our conventions in (\ref{2.5})
we have that $\bar{e}=\bar{e}_k$ and $\dim_\k{\mathcal
O}_\k(\bar{e})=2d(e)$. By well-known results of Spaltenstein and
Steinberg, all irreducible components of the quasiaffine variety
${\mathcal O}_\k(\bar{e}_i)\cap\n_\k^+$ have dimension equal to
$(\dim {\mathcal O}_\k(\bar{e}_i))/2$; see [\cite{Ja1}, 10.6] for
example.
\begin{theorem}\label{thm1}
Under the above assumptions on $A$ the variety $V_p(\mu)$ contains
an irreducible component of maximal dimension which coincides with
the Zariski closure of an irreducible component of ${\mathcal
O}_\k(\bar{e})\cap \n_\k^+$.
\end{theorem}
\begin{proof}
Given a nonempty locally closed subset $Z$ of the affine space
$\g_\k$ we denote by $\dim_\k Z$ the maximal dimension of the
irreducible components of $Z$. It is well-known that $\dim_\k
V_p(\mu)=\dim_\k{\mathbb K}(V_p(\mu))$; see [\cite{Kr}, II, 4.2].
In conjunction with (\ref{VA-VP}) and (\ref{Krull}) this gives
$\dim_\k V_p(\mu)\ge d(e)$.

\smallskip

\noindent (a) Let $X$ be an arbitrary component of maximal dimension
of $V_p(\mu)$. Since $\n_\k^+=\,\bigcup_{1\le i\le t}{\mathcal
O}_\k(\bar{e}_i)\cap\n_\k^+$, there exist $j\in\{1,\ldots, t\}$ and
an irreducible component $C$ of ${\mathcal
O}_\k(\bar{e}_j)\cap\n_\k^+$ such that $C\cap X$ is Zariski dense in
$X$. Then $X\subseteq \overline{C}$, implying
\begin{eqnarray}\label{ineq1}
\quad\qquad\dim_\k {\mathcal
O}_\k(\bar{e}_j)=2\dim_\k\overline{C}\ge 2\dim_\k X=2\dim_\k
V_p(\mu)\ge 2d(e)=\dim_\k{\mathcal O}_\k(\bar{e}).
\end{eqnarray}
Choose $x\in X\cap{\mathcal O}_k(\bar{e}_j)$ and set
$\xi:=(x,\cdot\,)$. Since $x\in V_p(\mu)$, the
$U_\xi(\g_\k)$-module $L_p^\xi(\mu)$ is nonzero; see
Lemma~\ref{lem1}. By the Kac--Weisfeiler conjecture proved in
[\cite{P95}], we have
$$p^{(\dim_\k {\mathcal O}_\k(\bar{e}_j))/2}\mid\dim_\k
L_p^\xi(\mu),$$ whilst Lemma~\ref{lem2} says that $\dim_\k
L_p^\xi(\mu)\le Dp^{d(e)}$ for some $D<p$. This enables us to
deduce that $\dim_\k{\mathcal O}_\k(\bar{e}_j)\le 2d(e).$ In
conjunction with (\ref{ineq1}) this implies that $\overline{C}=X$
and $\dim_\k X=\dim_\k V_p(\mu)=d(e)$ (one should keep in mind
that $\overline{C}$ is irreducible).

If $y\in {\mathcal O}_\k(\bar{e}_j)$ and $\lambda\in\k^{\times},$
then $\lambda y\in{\mathcal O}_\k(\bar{e}_j)$. From this it follows
that the Zariski closure of ${\mathcal O}_\k(\bar{e}_j)\cap\n_\k^+$
is a conical variety. But then so is $\overline{C}$, forcing
${\mathbb K}(X)=X$. Thus, {\it all} irreducible components of
maximal dimension of $V_p(\mu)$ are conical.

\smallskip

\noindent (b) Let $X_1,\ldots, X_c$ be the $d(e)$-dimensional
irreducible components of $V_p(\mu)$. Then $V_p(\mu)=X_1\cup\ldots
\cup X_c\cup Y$, where $Y$ is a closed subset of $V_p(\mu)$ with
$\dim_\k Y<d(e)$. It follows from [\cite{Kr}, II, 4.2] and part~(a)
of this proof that $${\mathbb K}(V_p(\mu))\,=\,X_1\cup\ldots\cup
X_c\cup{\mathbb K}(Y),\qquad \dim_\k{\mathbb K}(Y)<d(e).$$ In view
of (\ref{Krull}) and (\ref{VA-VP}) this shows that every irreducible
component of maximal dimension of ${\mathcal V}_\g L_p(\mu)$ is of
the form $X_i$ for some $i\le c$.

\smallskip

\noindent (c) Now suppose that the component $X=\overline{C}$ from
part~(a) of this proof is also an irreducible component of
${\mathcal V}_\g L_p(\mu)$. Such a component exists by our
discussion in part~(b). Recall from (\ref{2.4}) the set
$\{\psi_i\,|\,\,i\in I\}\subset\gr\,{\mathcal I}\cap S(\g_A)$ and
put $\bar{\psi}_i:=\psi_i\otimes 1\in S(\g_\k)$. Since
$\psi_i\in{\rm Ann}_{S(\g_A)}\,\gr\,L_A(\mu)$, we have
$\bar{\psi}_i\in {\rm Ann}_{S(\g_\k)}\,(L_A\otimes_A\k)={\mathcal
J}_p$. Hence
\begin{eqnarray}\label{psi}
{\mathcal V}_\g L_p(\mu)\,\subseteq \,\textstyle{\bigcap}_{i\in
I}\,V(\bar{\psi}_i).
\end{eqnarray}
Let $\Psi_A$ denote the $A$-span of the $\psi_i$'s in $S(\g_A)$,
and $\Psi_\k$ the $\k$-span of the $\bar{\psi}_i$'s in $S(\g_\k)$.
By our choices in (\ref{2.4}), the $A$-module $\Psi_A$ is stable
under the adjoint action of the hyperalgebra $U_\Z$. From this it
follows that the subspace $\Psi_\k\subset S(\g_\k)$ is invariant
under the adjoint action of the algebraic group $G_\k$.

Recall that $\overline{C}$ is a component of ${\mathcal
O}_\k(\bar{e}_j)\cap\n_\k^+$ and $\dim_\k  {\mathcal
O}_\k(\bar{e}_j)=2d(e)=2d(e_k)$. By our assumptions in
(\ref{2.5}), $\dim_{\mathbb C} {\mathcal O}(e_j)=2d(e)$, whilst
(\ref{psi}) shows that the $\bar{\psi}_i$'s vanish on
$\overline{C}$. Since $\Psi_\k={\rm span}\{\bar{\psi}_i\,|\,\,i\in
I\}$ is $G_\k$-stable, all $\bar{\psi}_i$ vanish on $(\Ad
G_\k)\cdot\overline{C}$. Therefore, $\bar{\psi}_i(\bar{e}_j)=0$
for all $i\in I$. But then $\psi_i(e_j)=0$ for all $i\in I$
forcing ${\mathcal O}(e_j)\subseteq \overline{{\mathcal O}(e_k)}$;
see (\ref{2.5}). Comparing dimensions now yields ${\mathcal
O}(e_j)={\mathcal O}(e_k)$, hence ${\mathcal
O}_\k(\bar{e}_j)={\mathcal O}_\k(\bar{e}_k)$. Thus, $X$ contains a
component of ${\mathcal O}_\k(\bar{e}_k)\cap\n_\k^+$, completing
the proof.
\end{proof}
\subsection{}\label{3.6} Recall that the group $G_\k$ acts on $U(\g_\k)$ as
algebra automorphisms (this action extends the adjoint action of
$G_\k$ on $\g_\k$). It is well-known (and easily seen) that for
every $\eta\in \g_\k^*$ and $g\in G_\k$ we have
$g(I_\eta)=I_{g\cdot\eta}$, where $g\cdot\eta=({\rm Ad}^* g)\eta$.
Thus, $g$ induces an algebra isomorphism
$U_\eta(\g_\k)\,\stackrel{\sim}{\longrightarrow}\,U_{g\cdot\eta}(\g_\k)$.
We set $\bar{\chi}=(\bar{e},\cdot\,)$.
\begin{corollary}\label{cor1}
Under the above assumptions on $A$, there exists a positive integer
$D=D(\mu)$ such that:
\begin{itemize}
\item[(i)] $D<p$ for all $p\in\pi(A)$;

\smallskip
\item[(ii)] for every $p\in\pi(A)$ the reduced enveloping algebra
$U_{\bar{\chi}}(\g_\k)$ has a simple module of dimension
$kp^{d(e)}$, where $k\le D$.
\end{itemize}
\end{corollary}
\begin{proof}
It follows from Theorem~\ref{thm1} that $g\cdot\bar{\chi}\in
V_p(\mu)$ for some $g\in G_\k$. Combining our earlier remark with
Lemmas~\ref{lem1} and \ref{lem2} we see that
$U_{\bar{\chi}}(\g_\k)$ has a nonzero module of dimension $\le
Dp^{d(e)}$, call it $E$. Since $D<p$, the Kac--Weisfeiler
conjecture (confirmed in [\cite{P95}]) shows that each composition
factor $L$ of the $\g_\k$-module $E$ has dimension $kp^{d(e)}$ for
some $k=k(L)\le D$.
\end{proof}
\section{\bf Non-restricted representations and finite $W$-algebras}
\subsection{}\label{4.1}
Let $e\in\g_\Z$ be as before, and  choose $f,h\in\g_{\mathbb Q}$
such that $(e,h,f)$ is an $\sl_2$-triple in $\g_{\mathbb Q}$. Let
$\g(i)=\{x\in\g\,|\,\,[h,x]=ix\}$, a subspace of $\g$ defined over
$\mathbb Q$. Note that $e\in\g(2)$, $f\in\g(-2)$, and
$\g=\bigoplus_{i\in\,\Z}\,\g(i)$. We define a skew-symmetric
bilinear form $\la\,\cdot\,,\,\cdot\,\ra$ on $\g(-1)$ by setting
$\la x,y\ra:=(e,[x,y])$ for all $x,y\in\g(-2)$. By the
$\sl_2$-theory, this skew-symmetric form is nondegenerate; hence
$\dim \g(-1)=2s$ for some $s\in \Z_+$.

There is a basis $B=\{z_1',\ldots,z_s',z_1,\ldots,z_s\}$ of $\g(-1)$
contained in $\g_\mathbb Q$ and such that
$$\la z_i',z_j\ra=\delta_{ij},\qquad\ \la z_i,z_j\ra\,=\,\la z_i',z_j'\ra=0
\qquad\quad\ (1\le i,j\le s).$$ Given an admissible subring $A$ of
$\mathbb Q$ we set $\g_A(i):=\g_A\cap\g(i)$, a finitely generated,
torsion-free $A$-module. {\it We shall assume, in addition, that
$\g_A=\bigoplus_{i\in\,\Z}\,\g_A(i)$, that each $\g_A(i)$ is a
freely generated over $A$ by a basis of the vector space $\g(i)$,
and that $B$ is a free basis of the $A$-module $\g_A(-1)$}. It is
straightforward to see that admissible subalgebras $A$ satisfying
these conditions exist.

Let $\g(-1)^0$ denote the $\mathbb C$-span of $z_1',\ldots, z_s'$
and put $\m:=\g(-1)^0\oplus\sum_{i\le  -2}\,\g(i)$. Note that $\m$
is a nilpotent Lie subalgebra of dimension $d(e)$ in $\g$ and
$\chi$ vanishes on the derived subalgebra of $\m$; see
[\cite{P02}] for more detail. Put $\m_A:=\g_A\cap\m$. It follows
from the above assumptions that $\m_A$ is a free $A$-module and a
direct summand of $\g_A$. More precisely,
$\m_A=\g_A(-1)^0\oplus\sum_{i\le -2}\,\g_A(i)$ where
$\g_A(-1)^0=\g_A\cap\g(-1)=Az_1'\oplus\cdots\oplus Az_s'$.

Let $\g_e$ denote the centralizer of $e$ in $\g$. Clearly, $\g_e$
coincides with the stabiliser
$\z_\chi=\{x\in\g\,|\,\,\chi([x,\g])=0\}$. By the $\sl_2$-theory,
the map $\ad e\colon\,\g_{\mathbb Q}(i)\rightarrow \g_{\mathbb
Q}(i+2)$ is surjective for all $i\ge 0$ and $\g_e=\bigoplus_{i\ge
0}\,\z_\chi(i)$ where $\z_\chi(i)=\z_\chi\cap\g(i)$. Enlarging $A$
if need be {\it we shall assume that $\g_A(i+2)=[e,\g_A(i)]$ for
all $i\ge 0$}. Then there is a basis $x_1,\ldots,
x_r,x_{r+1},\ldots, x_m$ of the free $A$-module $\bigoplus_{i\ge
0}\,\g_A(i)$ such that $x_i\in\g_A(n_i)$ for some $n_i\in\Z_+$,
where $1\le i\le m$, and $x_1,\ldots, x_r$ is a free basis of the
$A$-module $\g_A\cap\z_\chi$; see [\cite{P02}, 4.2 \& 4.3] for
more detail.
\subsection{}\label{4.2}
Let ${\mathbb C}_\chi={\mathbb C}1_\chi$ be a $1$-dimensional
$\m$-module such that $x\cdot 1_\chi=\chi(x)1_\chi$ form all
$x\in\m$, and $A_\chi=A1_\chi$. Given $({\bf a},{\bf
b})\in\Z_+^m\times\Z_+^s$ we denote by $x^{\bf a}z^{\bf b}$ the
monomial $x_1^{a_1}\cdots x_m^{a_m}z_1^{b_1}\cdots z_s^{b_s}$ in
$U(\g)$. Set $Q_\chi:=U(\g)\otimes_{U(\m)}{\mathbb C}_\chi$ and
$Q_{\chi,\,A}:=U(\g_A)\otimes_{U(\m_A)}A_\chi$. From our
discussion in (\ref{4.1}) it follows that $Q_{\chi,\,A}$ is an
$A$-lattice in $Q_\chi$ with free basis consisting of all $x^{\bf
i}z^{\bf j}\otimes 1_\chi$ with $({\bf i},{\bf
j})\in\Z_+^m\times\Z_+^s$. Moreover, $Q_{\chi,\,A}$ is invariant
under the action of $\g_A$. Given $({\bf a},{\bf
b})\in\Z_+^m\times Z_+^s$ we set
$$|({\bf a},{\bf b})|_e\,:=\,\,\sum_{i=1}^m\,a_i(n_i+2)+\sum_{i=1}^s\,b_i.$$
According to [\cite{P02}, Theorem~4.6], the algebra
$H_\chi:=(\End_\g\,Q_\chi)^{\rm op}$ is generated over $\mathbb C$
by endomorphisms $\Theta_1,\ldots,\Theta_r$ such that
\begin{eqnarray}\label{lam} \qquad\ \,\Theta_k(1_\chi)\,=\,\Big(x_k+\sum_{0<|({\bf i},{\bf j})|_e\le
n_k+2}\,\lambda_{{\bf i},\,{\bf j}}^k\,x^{\bf i}z^{\bf
j}\Big)\otimes 1_\chi,\qquad\quad \ 1\le k\le r,\end{eqnarray} where
$\lambda_{{\bf i},\,{\bf j}}^k\in\mathbb Q$ and $\lambda_{{\bf
i},\,{\bf j}}^k=0$ if either $|({\bf i},\,{\bf j})|_e=n_k+2$ and
$|{\bf i}|+|{\bf j}|=1$ or ${\bf i}\ne {\bf 0}$, ${\bf j}={\bf 0}$,
and $i_l=0$ for $l>r$.

The monomials $\Theta_1^{i_1}\cdots\Theta_r^{i_r}$ with
$(i_1,\ldots, i_r)\in\Z_+^r$ form a PBW basis of the vector space
$H_\chi$. The monomial $\Theta_1^{i_1}\cdots\Theta_r^{i_r}$ is
said to have {\it Kazhdan degree} $\sum_{i=1}^r\,a_i(n_i+2)$. For
$k\in\Z_+$ we let $H_\chi^k$ denote the $\mathbb C$-span of all
monomials $\Theta_1^{i_1}\cdots\Theta_r^{i_r}$ of Kazhdan degree
$\le k$. Then $(H_\chi^k)_{k\ge 0}$ is an increasing filtration of
the algebra $H_\chi$, called the {\it Kazhdan filtration} of
$H_\chi$. The corresponding graded algebra $\gr\,H_\chi$ is a
polynomial algebra in $\gr\,\Theta_1,\ldots,\gr\,\Theta_r$. It is
immediate from [\cite{P02}, 4.5] that there exist polynomials
$F_{ij}\in{\mathbb Q}[X_1,\ldots, X_r]$, where $1\le i<j\le r$,
such that
\begin{eqnarray}\label{relations}
\qquad\
[\Theta_i,\Theta_j]\,=\,F_{ij}(\Theta_1,\ldots,\Theta_r)\qquad\quad\
\ (1\le i<j\le r).
\end{eqnarray}
Moreover, if $[x_i,x_j]\,\,=\,\,\sum_{k=1}^r \alpha_{ij}^k\, x_k$ in
$\z_\chi$, then
$$F_{ij}(\Theta_1,\ldots,\Theta_r)\,\equiv\,\sum_{k=1}^r\alpha_{ij}^k\Theta_k
+q_{ij}(\Theta_1,\ldots,\Theta_r)\ \ \ \big({\rm
mod}\,H_\chi^{n_i+n_j}\big)$$ where the initial form of $q_{ij}\in
{\mathbb Q}[X_1,\ldots, X_r]$ has total degree $\ge 2$ whenever
$q_{ij}\ne 0$.
\begin{lemma}\label{lem3}
The algebra $H_\chi$ is generated by $\Theta_1,\ldots,\Theta_r$
subject to the relations (\ref{relations}).
\end{lemma}
\begin{proof}
Let $\bf I$ be the two-sided ideal of the free associative algebra
${\mathbb C}\la X_1,\ldots,X_r\ra$ generated by all
$[X_i,X_j]-F_{ij}(X_1,\ldots, X_r)$ with $1\le i<j\le r$. Let
$\hat{H}_\chi:={\mathbb C}\la X_1,\ldots,X_r\ra/{\bf I}$, and let
$\hat{\Theta}_i$ be the image of $X_i$ in $\hat{H}_\chi$. There is
a natural algebra epimorphism $\hat{H}_\chi\twoheadrightarrow
H_\chi$ sending $\hat{\Theta}_i$ to $\Theta_i$ for all $i$; we
call it $\psi$. For $k\in\Z_+$ let $\hat{H}_\chi^k$ denote the
$\mathbb C$-span of all products $\hat{\Theta}_{j_1}\cdots
\hat{\Theta}_{j_l}$ with $(n_{j_1}+2)+\cdots+(n_{j_l}+2)\le k$.
Let $\hat{H}_\chi'$ denote the $\mathbb C$-span of all monomials
$\hat{\Theta}_1^{i_1}\cdots\hat{\Theta}_r^{i_r}$ with
$(i_1,\ldots, i_r)\in\Z_+^r$. Taking into account the relations
(\ref{relations}) and arguing by upward induction on $k$ and
downward induction on $l$, the number of factors in
$\hat{\Theta}_{j_1}\cdots\hat{\Theta}_{j_l}\in\hat{H}_\chi^k$, one
observes that $\hat{H}_\chi=\hat{H}_\chi'$. Since the monomials
$\Theta_1^{i_1}\cdots\Theta_r^{i_r}$ with
$(i_1,\ldots,i_r)\in\Z_+^r$ are linearly independent in $H_\chi$,
this implies that $\psi$ is injective. But then $\hat{H}_\chi\cong
H_\chi$, and our proof is complete.
\end{proof}
Enlarging $A$ if necessary {\it we shall assume that all
$\lambda_{{\bf i},{\bf j}}^k$ in (\ref{lam}) and all coefficients
of the polynomials $F_{ij}$ in (\ref{relations}) are in $A$}.
\subsection{}\label{4.3}
Let $H_{\chi,\,A}$ denote the $A$-span of all monomials
$\Theta_1^{i_1}\cdots\Theta_r^{i_r}$ with
$(i_1,\ldots,i_r)\in\Z_+^r$. Our assumptions on $A$ in (\ref{4.2})
show that $H_{\chi,\,A}$ is an $A$-subalgebra of $H_\chi$
contained in $(\End_{\g_A}\,Q_{\chi,\,A})^{\rm op}$. We set
$\m_{\k}:=\m_A\otimes_A\k$ and
$Q_{\chi,\,\k}:=U(\g_\k)\otimes_{U(\m_\k)}\k_{\bar{\chi}}$ where
$\k_{\bar{\chi}}=A_\chi\otimes_A\k\,=\,\k1_{\bar{\chi}}$. It
follows from our assumptions in (\ref{4.1}) that
$\bar{x}_1,\ldots, \bar{x}_r$ form a basis of the stabiliser
$\z_{\bar{\chi}}$ of $\bar{\chi}$ in $\g_\k$ and that $\m_{\k}$ is
a nilpotent subalgebra of dimension $d(e)$ in $\g_\k$. Also,
$\k1_{\bar{\chi}}$ is a $1$-dimensional $\m_\k$-module with the
property that $x(1_{\bar{\chi}})=\bar{\chi}(x)1_{\bar{\chi}}$ for
all $x\in\m_\k$. Let $I_{\bar{\chi}}$ be the two-sided ideal of
$U(\g_\k)$ generated by all $x^p-x^{[p]}-\bar{\chi}(x)^p1$ with
$x\in\g_\k$, and set
$Q_{\bar{\chi}}^{[p]}:=Q_{\chi,\,\k}/I_{\bar{\chi}}Q_{\chi,\,\k}$,
a $\g_\k$-module with $p$-character $\bar{\chi}$. Each
$\g_\k$-endomorphism $\Theta_i\otimes 1$ of $Q_{\chi,\,\k}=
Q_{\chi,\,A}\otimes_{A}\k$ preserves the submodule
$I_{\bar{\chi}}Q_{\chi,\,\k}$, hence induces a
$\g_\k$-endomorphism of $Q_{\bar{\chi}}^{[p]}$; call it
$\theta_i$. We denote by $H_{\bar{\chi}}^{[p]}$ the associative
$\k$-algebra $\big(\!\End_{\g_\k}\,Q_{\bar{\chi}}^{[p]}\big)^{\rm
op}$.

Note that the $\g_A$-module $Q_{\chi,\,A}$ is isomorphic to
$U(\g_A)/U(\g_A)\n_{\chi,\,A}$, where $\n_{\chi,\,A}$ is the
$A$-span in $U(\g_A)$ of all $x-\chi(x)$ with $x\in\m_{A}$. It
follows that $H_{\chi,\,A}$ identifies with an $A$-subalgebra of
$\big(U(\g_A)/U(\g_A)\n_{\chi,\,A}\big)^{\ad\m_A}$. According to
[\cite{P02}, 2.3] we have that
$H_{\bar{\chi}}^{[p]}\,\cong\,\big(U_{\bar{\chi}}(\g_\k)/
U_{\bar{\chi}}(\g_\k)\n_{\chi,\,\k}\big)^{\ad\m_\k}$ where
$\n_{\chi,\,\k}=\n_{\chi,\,A}\otimes_A\k$.
\begin{prop}\label{prop1}
The following are true:
\begin{itemize}
\item[(i)] $Q_{\bar{\chi}}^{[p]}\,\cong\,
U_{\bar{\chi}}(\g_\k)\otimes_{U_{\bar{\chi}}(\m_\k)}\k_{\bar{\chi}}$
as $\g_\k$-modules;

\smallskip
\item[(ii)] $Q_{\bar{\chi}}^{[p]}$ is a projective generator for
$U_{\bar{\chi}}(\g_\k)$ and
$U_{\bar{\chi}}(\g_\k)\,\cong\,\Mat_{p^{d(e)}}\big(H_{\bar{\chi}}^{[p]}\big)$;

\smallskip

\item[(iii)] the monomials $\theta_1^{i_1}\cdots\theta_r^{i_r}$ with
$0\le i_k\le p-1$ form a $\k$-basis of $H_{\bar{\chi}}^{[p]}$.
\end{itemize}
\end{prop}
\begin{proof}
Let $\bar{1}_{\bar{\chi}}$ denote the image of $1_{\bar{\chi}}$ in
$Q_{\bar{\chi}}^{[p]}$. By the universality property of induced
modules, there is a surjective homomorphism
$\tilde{\alpha}\colon\,Q_{\chi,\,\k}=U(\g_\k)\otimes_{U(\m_\k)}\k_{\bar{\chi}}\twoheadrightarrow
U_{\bar{\chi}}(\g_\k)\otimes_{U_{\bar{\chi}}(\m_\k)}\k_{\bar{\chi}}$.
Since $I_{\bar{\chi}}Q_{\chi,\,\k}\subseteq \ker\tilde{\alpha}$,
it gives rise to an epimorphism
$\alpha\colon\,Q_{\bar{\chi}}^{[p]}\twoheadrightarrow
U_{\bar{\chi}}(\g_\k)\otimes_{U_{\bar{\chi}}(\m_\k)}\k_{\bar{\chi}}$.
On the other hand, the $U_{\bar{\chi}}(\g_\k)$-module
$Q_{\bar{\chi}}^{[p]}$ is generated by its $1$-dimensional
$U_{\bar{\chi}}(\m_\k)$-submodule
$\k\bar{1}_{\bar{\chi}}\cong\k_{\bar{\chi}}$. The universality
property of induced $U_{\bar{\chi}}(\g_\k)$-modules now yields
that there is a surjection
$\beta\colon\,U_{\bar{\chi}}(\g_\k)\otimes_{U_{\bar{\chi}}(\m_\k)}\k_{\bar{\chi}}\twoheadrightarrow
Q_{\bar{\chi}}^{[p]}.$ Comparing dimensions shows that $\alpha$ is
an isomorphism, proving (i). As $\m_\k\cap\z_{\bar{\chi}}=0$ by
our assumptions in (\ref{4.1}), $\m_\k$ is a
$\bar{\chi}$-admissible subalgebra of dimension $d(e)$ in $\g_\k$;
see [\cite{P02}, 2.3 \& 2.6]. Now (ii) and (iii) follow from
[\cite{P02}, Theorems~2.3 \& 3.4].
\end{proof}
\subsection{}\label{4.4} As $\dim\g(i)=\dim\g(-i)$ for all $i$
and $r=\dim\g(0)+2s$ by the $\sl_2$-theory, we have that
$\dim\,\bigoplus_{i\le -2}\g(i)=m-r$ and $d(e)=m-r+s$. Set
$$X_i=\left\{
\begin{array}{ll}
z_i&\mbox{if $\ 1\le i\le s$},\\
x_{r-s+i}&\mbox{if $\ s+1\le i\le m-r+s$},
\end{array}\right.
$$ and let
$\bar{X}_i$ denote the image of $X_i$ in $\g_\k=\g_A\otimes_A\k$.
Since $H_{\chi,\,A}\subseteq(\End_{\g_A}Q_{\chi,A})^{\rm op}$, we
can regard $Q_{\chi,\,A}$ as a right $H_{\chi,\,A}$-module, whilst
from the definition of $H_{\bar{\chi}}^{[p]}$ it is clear that
$Q_{\bar{\chi}}^{[p]}$ is a right $H_{\bar{\chi}}^{[p]}$-module.
Proposition~\ref{prop1} shows that these two module structures are
compatible, that is the latter can be obtained obtained from the
former by reducing $Q_{\chi,\,A}$ modulo $p$ and then by passing
from $Q_{\chi,\,\k}$ to its quotient $Q_{\bar{\chi}}^{[p]}$.

Given $k\ge 0$ we let $Q_{\chi,\,A}^k$ denote the $A$-submodule of
$Q_{\chi,\,A}$ spanned by all $x^{\bf i}z^{\bf j}\otimes 1_\chi$
with $|({\bf i},{\bf j})|_e\le k$. For ${\bf a}\in \Z_+^{d(e)}$ we
define $X^{\bf a}:=X_1^{a_1}\cdots X_{d(e)}^{a_{d(e)}}\otimes
1_\chi$ and $\bar{X}^{\bf a}:=\bar{X}_1^{a_1}\cdots
\bar{X}_{d(e)}^{a_{d(e)}}\otimes \bar{1}_{\bar{\chi}}$, elements
of $Q_{\chi,\,A}$ and $Q_{\bar{\chi}}^{[p]}$, respectively.
\begin{lemma}\label{lem4} The right modules
$Q_{\chi,\,A}$ and $Q_{\bar{\chi}}^{[p]}$ are free over
$H_{\chi,\,A}$ and $H_{\bar{\chi}}^{[p]}$, respectively. More
precisely,
\begin{itemize}
\item[(i)] the set $\big\{X^{\bf a}\otimes 1_\chi\,|\,\,{\bf
a}\in\Z_+^{d(e)}\big\}$ is a free basis of the
$H_{\chi,\,A}$-module $Q_{\chi,\,A}$;
\smallskip

\item[(ii)] the set $\big\{\bar{X}^{\bf a}\otimes
\bar{1}_{\bar{\chi}}\,|\,\, 0\le a_i\le p-1\big\}$ is a free basis
of the $H_{\bar{\chi}}^{[p]}$-module $Q_{\bar{\chi}}^{[p]}$.
\end{itemize}
\end{lemma}
\begin{proof}
For ${\bf a}\in\Z_+^r$ set $\Theta^{\bf a}:= \Theta_1^{a_1}\cdots
\Theta_r^{a_r}$. Using (\ref{lam}) and induction on the Kazhdan
degree $k=\sum_{i=1}^ra_i(n_i+2)$ of $\Theta^{\bf a}$ it is easy to
observe that
$$\Theta^{\bf a}(1_\chi)\,\equiv\, x_1^{a_1}\cdots
x_r^{a_r}\otimes 1_\chi+\ \sum_{|({\bf i},{\bf j})|_e=k,\ |{\bf
i}|+|{\bf j}|>|\bf a|}\,\gamma_{{\bf i},{\bf j}}\,x^{\bf i}z^{\bf
j}\otimes 1_\chi\quad\ (\mathrm{mod}\ Q_\chi^{k-1})$$ for some
$\gamma_{{\bf i},{\bf j}}\in A$; cf. [\cite{P02}. p.~27]. Using
this relation and arguing by double induction on $|({\bf i},{\bf
j})|_e$ and $|{\bf i}|+|{\bf j}|$ (upward on $|({\bf i},{\bf
j})|_e$ and downward on $|{\bf i}|+|{\bf j}|$) one observes that
every $x^{\bf i}z^{\bf j}\otimes 1_\chi$ belongs to the
$A$-submodule of $Q_{\chi,\,A}$ spanned by the vectors $X^{\bf
a}\Theta^{\bf b}(1_\chi)$ with ${\bf a}\in\Z_+^{d(e)}$ and ${\bf
b}\in\Z_+^r$. Since the latter vectors are linearly independent in
$H_\chi$ in view of [\cite{Sk}, pp.~52, 53], statement~(i)
follows. Statement~(ii) is proved  similarly; see [\cite{P02},
3.4] for more detail.
\end{proof}
\subsection{}\label{4.5} Let $\{u_i\,|\,\,i\in I\}$ be the
finite generating set of the primitive ideal $\mathcal I$
discussed in (\ref{2.4}). Recall that $u_i\in U(\g_A)$ for all
$i\in I$ and the $A$-span of the $u_i$'s is invariant under the
adjoint action of $\g_A$ on $U(\g_A)$. For $i\in I$ we denote by
$\bar{u}_i$ the image of $u_i$ in $U(\g_\k)=U(\g_A)\otimes_A\k$.
By construction, the $\k$-span of the $\bar{u}_i$'s is invariant
under the adjoint action of $\g_\k$ on $U(\g_\k)$. Let
$\varphi_\chi\colon\, U(\g_A)\twoheadrightarrow
Q_{\chi,\,A}=U(\g_A)/U(\g_A)\n_{\chi,\,A}$ denote the canonical
homomorphism, and $\bar{\varphi}_\chi$ the induced epimorphism
from $U_{\bar{\chi}}(\g_\k)$ onto $Q_{\bar{\chi}}^{[p]}$; see
(\ref{4.3}) for more detail.  It follows from Lemma~\ref{lem4}(i)
that there exists a finite subset $C$ of $\Z_+^{d(e)}$ such that
\begin{eqnarray}\label{u-h}\varphi_\chi(u_i)\,=\, \sum_{{\bf c}\in\, C}\,X^{\bf c}
h_{i,\,{\bf c}}(1_\chi)\qquad\qquad \ (h_{i,\,\bf c}\in
H_{\chi,\,A},\ \,i\in I).
\end{eqnarray}
It follows from Proposition~\ref{prop1}(iii) that the $\k$-algebra
$H_{\bar{\chi}}^{[p]}$ is a homomorphic image of the $\k$-algebra
$H_{\chi,\, \k}:=H_{\chi,\,A}\otimes_A\k$.  Let $\bar{h}_{i,\,\bf
c}$ denote the image of $h_{i,\,\bf c}\otimes 1$ in
$H_{\bar{\chi}}^{[p]}$. Then it follows from (\ref{u-h}) that
\begin{eqnarray}\label{u-h-p}\bar{\varphi}_\chi(\bar{u}_i)\,=\, \sum_{{\bf c}\in\, C}\,
\bar{X}^{\bf c} \bar{h}_{i,\,{\bf
c}}(\bar{1}_{\bar{\chi}})\qquad\qquad \ (\forall\,i\in I).
\end{eqnarray}
Let $c:=\max_{{\bf c}\in\, C}|{\bf c}|$. We now impose our final
assumption on $A$: {\it we assume that $c!$ is invertible in $A$}.
This assumption ensures that the components of all tuples in $C$
are smaller that any prime in $\pi(A)$.
\begin{prop}\label{prop2}
Under the above assumption on $A$, for every $p\in\pi(A)$ there
exists a positive integer $k=k(p)\le D=D(\mu)$ such that the
associative algebra $H_{\bar{\chi}}^{[p]}$ has an irreducible
$k$-dimensional representation $\rho$ with the property that
$\rho(\bar{h}_{i,\,{\bf c}})=0$ for all ${\bf c}\in C$ and all
$i\in I$.
\end{prop}
\begin{proof}
Let $p\in\pi(A)$ and $\k=\overline{\mathbb F}_p$. By
Lemma~\ref{lem2}, Theorem~\ref{thm1} and Corollary~\ref{cor1},
there exists $g\in G_\k$ such that any composition factor $V$ of
the $\g_\k$-module $L_p^{g\cdot\bar{\chi}}(\mu)$ has dimension
$kp^{d(e)}$ for some $k=k(V)\le D$. Since $u_i\in{\rm
Ann}_{U(\g_A)}\,L_A(\mu)$ for all $i\in I$, the elements
$\bar{u}_i\in U(\g_\k)$ annihilate $L_p(\mu)=L_A(\mu)\otimes_A\k$.
Consequently, all $\bar{u}_i$ annihilate
$L_p^{g\cdot\bar{\chi}}(\mu)=L_p(\mu)/I_{g\cdot\bar{\chi}}L_p(\mu)$,
and hence $V$. Let $V'=\{v'\,|\,\,\,v\in V\}$ be another copy of
the vector space $V$. We give $V'$ a $\g_\k$-module structure by
setting $x.v'\,:=\,((\Ad g)^{-1}x . v)'$ for all $x\in\g_\k$ and
all $v'\in V'.$ Since all elements $((\Ad g)x)^p-((\Ad
g)x)^{[p]}-\chi(x)^p1$ annihilate $V$, the $\g_\k$-module $V'$ has
$p$-character $\bar{\chi}$. By construction, all elements $(\Ad
g)\bar{u}_i$ annihilate $V'$.

Recall that the $\Z$-span of $\{u_i\,|\,\,i\in I\}$ is invariant
under the adjoint action of the hyperalgebra $U_\Z$ on $U(\g_\Z)$;
see (\ref{2.4}). This implies that the $\k$-span of the
$\bar{u}_i$'s is invariant under the adjoint action of $G_\k$ on
$U(\g_\k)$. In conjunction with our preceding remark this implies
that $\bar{u}_i\in{\rm Ann}_{U(\g_\k)}\,V'$ for all $i\in I$. Let
$$V'_0\,=\,\{v'\in V'\,|\,\,\,x.v'=\,\bar{\chi}(x)v'\ \, \mbox{for all
}\,x\in\m_\k\},$$ the subspace of $\m_\k$-Whittaker vectors in
$V'$. Since $H_{\bar{\chi}}^{[p]}\,\cong
\big(U_{\bar{\chi}}(\g_\k)/U_{\bar{\chi}}(\g_\k)\n_{\chi,\,\k}\big)^{\ad\m_\k}$,
the algebra $H_{\bar{\chi}}^{[p]}$ acts on $V'_0$. Since $\m_\k$
is a $\bar{\chi}$-admissible subalgebra of dimension $d(e)$ in
$\g_\k$, it follows from [\cite{P02}, Theorem~2.4] that $V'_0$ is
an irreducible $k$-dimensional $H_{\bar{\chi}}^{[p]}$-module. We
denote by $\rho$ the corresponding representation of
$H_{\bar{\chi}}^{[p]}$.

Let $v_1',\ldots, v_k'$ be a basis of $V_0'$ and
$V'':=Q_{\bar{\chi}}^{[p]}\otimes_{H_{\bar{\chi}}^{[p]}}V_0'$, a
$\g_\k$-module with $p$-character $\bar{\chi}$. It follows from
Lemma~\ref{lem4}(ii) that the vectors $\bar{X}^{\bf a}\otimes
v_j'$ with $0\le a_i\le p-1$ and $1\le j\le k$ form a basis of
$V''$ over $\k$. Since $V'$ is irreducible, there is a
$\g_\k$-module epimorphism $\tau\colon\,V''\twoheadrightarrow V'$
sending $v'\otimes 1$ to $v'$ for all $v'\in V_0'$. Since $\dim_\k
V'=kp^{d(e)}$, comparing dimensions yields that $\tau$ is an
isomorphism. Let $\tilde{\rho}$ denote the representation of
$U_{\bar{\chi}}(\g_\k)$ in $\End_\k V''$. As the left ideal
$U_{\bar{\chi}}(\g_\k)\n_{\chi,\,\k}$ of $U_{\bar{\chi}}(\g_\k)$
annihilates the subspace $V'_0\otimes 1$ of $V''$, it follows from
(\ref{u-h-p}) that
$$
0\,=\,\tilde{\rho}(\bar{u}_i)(v'\otimes
1)\,=\,\tilde{\rho}(\bar{\varphi}_\chi(\bar{u}_i))(v'\otimes 1)
\,=\,\,\sum_{{\bf c}\in\, C}\,\,\bar{X}^{\bf
c}\otimes\rho(\bar{h}_{i,\,{\bf c}})(v')
$$
for all $v'\in V'_0$. As the nonzero vectors of the form
$\bar{X}^{\bf c}\otimes\rho(\bar{h}_{i,\,{\bf c}})(v')$ with $v'$
fixed are linearly independent in $V''$  by our final assumption
on $A$, we now derive that $\rho(\bar{h}_{i,\,{\bf c}})=0$ for all
${\bf c}\in C$ and all $i\in I$. This completes the proof.
\end{proof}
\subsection{}\label{4.6}
Since $H_{\chi,\,A}$ is a free $A$-module with basis
$\{\Theta^{\bf i}\,|\,\,\,{\bf i}\in\Z_+^r\}$ there exist
polynomials $H_{i,\,{\bf c}}\in A[X_1,\ldots, X_r]$ such that
$h_{i,\,{\bf c}}=H_{i,\,{\bf c}}(\Theta_1,\ldots, \Theta_r)$ for
all ${\bf c}\in C$ and all $i\in I$. Let ${\mathcal I}_H$ denote
the two-sided ideal of $H_\chi$ generated by the $h_{i,\,{\bf
c}}$'s. In view of (\ref{relations}) and Lemma~\ref{lem3}, the
algebra $H_\chi/{\mathcal I}_H$ is isomorphic to the quotient of
the free associative algebra ${\mathbb C}\la X_1,\ldots, X_r\ra$
by its two-sided ideal generated by all elements
$[X_i,X_j]-F_{ij}(X_1,\ldots, X_r)$ with $1\le i<j\le r$ and all
elements $H_{{\bf c},\,l}(X_1,\ldots, X_r)$ with ${\bf c}\in C$
and $l\in I$. Given a natural number $k$ we let ${\mathcal Y}_k$
denote the set of all $r$-tuples $(M_1,\ldots, M_r)\in
\Mat_k({\mathbb C})^r$ satisfying the relations
\begin{eqnarray*}
[M_i,M_j]-F_{ij}(M_1,\ldots, M_r)&=&0\qquad\quad (1\le i<j\le r)\\
H_{{\bf c},\,l}(M_1,\ldots, M_r)&=&0 \qquad\quad ({\bf c}\in C,\,\,
l\in I)
\end{eqnarray*}
(for $F\in {\mathbb C}[X_1,\ldots,X_r]$, all monomials of
$F(M_1,\ldots,M_r)$ are understood to be taken with respect to the
matrix product in $\Mat_k(\mathbb C)$). The above discussion shows
that ${\mathcal Y}_k$ is the variety of all matrix representations
of degree $k$ of the algebra $H_\chi/{\mathcal I}_H$.

We regard ${\mathcal Y}_k$ as a Zariski closed subset of ${\mathbb
A}^{rk^2}({\mathbb C})$. More precisely, denote by ${\bf J}_k$ the
ideal of the polynomial algebra $P:={\mathbb
C}\,[x_{ab}^{(c)}\,|\,\,\,1\le a,b\le k,\, 1\le c\le r]$ generated
by the matrix coefficients of all $[M_i,M_j]-F_{ij}(M_1,\ldots,
M_r)$ and $H_{{\bf c},\,l}(M_1,\ldots, M_r)$, where $M_c$ is the
generic matrix $\big(x_{ab}^{(c)}\big)_{1\le a,b\le
k}\in\Mat_k(P)$ and $1\le c\le r$. Then ${\mathcal Y}_k$ is
nothing but the zero locus of ${\mathbf J}_k$ in ${\rm
Specm}\,P=\,{\mathbb A}^{rk^2}({\mathbb C})$.

Set $P_\Z:=\,\Z\,[x_{ab}^{(c)}\,|\,\,\,1\le a,b\le k,\, 1\le c\le
r]$. As all $F_{ij}$ and $H_{{\bf c},\,l}$ are in $A[X_1,\ldots,
X_r]$, the ideal ${\bf J}_k$ is generated by a finite set of
polynomials in $P_A:=P_{\Z}\otimes_{\Z}A$, say $f_1,\ldots, f_R$.
Given $p\in\pi(A)$ and $g\in P_A$ we write $^{p\!}{g}$ for the
image of $g\in P_A$ in $P_A\otimes_A\overline{\mathbb F}_p$, and
we denote by ${\mathcal Y}_k(\overline{\mathbb F}_p)$ the zero
locus of $^{p\!}f_1,\ldots, {^{p\!}f_R}$ in ${\mathbb
A}^{rk^2}({\overline{\mathbb F}_p})$.
\begin{prop}\label{prop3}
The algebra $H_\chi/{\mathcal I}_H$ has an irreducible
representation of dimension at most $D=D(\mu)$.
\end{prop}
\begin{proof}
In view of the above discussion it suffices to show that for some
$k\le D$ the variety ${\mathcal Y}_k$ has a point with
coefficients in $\overline{\mathbb Q}$, the algebraic closure of
$\mathbb Q$ in $\mathbb C$. Suppose this is not the case. Then
$g'_1f_1+\cdots+g'_Rf_R=1$ for some $g'_1,\ldots,g'_R\in
P_A\otimes_A \overline{\mathbb Q}$. Let $K\subset\mathbb C$ be a
finite Galois extension of $\mathbb Q$ containing all coefficients
of $g'_1,\ldots,g'_R$, and let $\Lambda$ be the ring of algebraic
integers of $K$. Rescaling  $g_1',\ldots, g_R'$ we find
$g_1,\ldots, g_R\in P_{\Lambda}:= P_{\Z}\otimes_{\Z}\Lambda$ and a
{\it nonzero} $\tilde{n}\in\Z$ such that $g_1f_1+\cdots+
g_Rf_R=\tilde{n}.$

It is well-known that $\Lambda$ is a Dedekind ring and the map ${\rm
Spec}\,\Lambda\rightarrow {\rm Spec}\,\Z$ induced by inclusion
$\Z\hookrightarrow\Lambda$ is surjective. For each $p\in\pi(A)$ we
choose ${\mathfrak P}\in{\rm Spec}\,\Lambda$ such that ${\mathfrak
P}\cap\Z=p\Z$ and denote by $\phi$ the composite
$$P_\Lambda\,\twoheadrightarrow\, P_\Lambda/{\mathfrak
P}P_\Lambda\,\hookrightarrow\, \overline{\mathbb
F}_p\,[x_{ab}^{(c)}\,|\,\,\,1\le a,b\le k,\, 1\le c\le
r]=P_A\otimes_A\overline{\mathbb F}_p.$$ Since
$\phi(f_i)={^{p\!}}f_i$ for all $i$ by our choice of $\mathfrak
P$, we have that
$$\phi(g_1)\cdot {{^p\!}}f_1+\cdots+ \phi(g_R)\cdot
{^{p\!}}f_R=\phi(\tilde{n}).$$ Since $\phi(\tilde{n})$ is nothing
but the image of $\tilde{n}$ in ${\mathbb F}_p$, this implies that
${\mathcal Y}_k(\overline{\mathbb F}_p)=\emptyset$ for all $k\le
D$ and almost all $p\in\pi(A)$. On the other hand, it follows from
Proposition~\ref{prop2} and the preceding discussion in
(\ref{4.5}) that that there is a positive integer $l\le D$ such
that ${\mathcal Y}_{l}(\overline{\mathbb F}_p)\ne \emptyset$ for
infinitely many $p\in\pi(A)$. Since $\tilde{n}$ has only finitely
many prime divisors, we reach a contradiction thereby completing
the proof.
\end{proof}
\subsection{}\label{4.7} We are now ready to prove the main
results of this note.
\begin{theorem}\label{main}
Let $\mathcal I$ be a primitive ideal of $U(\g)$ with rational
infinitesimal character and with $\mathcal{VA}({\mathcal
I})\,=\,\overline{\mathcal O}_\chi$. Then ${\mathcal I}\,=\,{\rm
Ann}_{U(\g)}(Q_\chi\otimes_{H_\chi} V)$ for some finite
dimensional irreducible $H_\chi$-module $V$.
\end{theorem}
\begin{proof}
By Proposition~\ref{prop3}, the algebra $H_\chi/{\mathcal I}_H$
has a finite- dimensional irreducible module, say $V$. We write
$\rho$ for the corresponding representation of $H_\chi$ in $\End
V$ and set $\widetilde{V}:=\,Q_\chi\otimes_{H_\chi}V$. Denote by
$\tilde{\rho}$ the representation of $U(\g)$ in
$\End\widetilde{V}$ and set $\tilde{\mathcal I}:= {\rm
Ann}_{U(\g)}\big(Q_\chi\otimes_{H_\chi} V\big)$. Skryabin's
theorem [\cite{Sk}] implies that $\widetilde{V}$ is an irreducible
$\g$-module, whilst [\cite{P07}, Theorem~3.1(ii)] says that
$\mathcal{VA}(\tilde{\mathcal I})=\overline{\mathcal{O}}_\chi$. In
view of (\ref{u-h}),
$$\tilde{\rho}(u_i)(1\otimes v)\,=\,\tilde{\rho}(\varphi_\chi(u_i))(1\otimes v)
\,=\, \sum_{{\bf c}\in\, C}\,X^{\bf c}\otimes \rho(h_{i,\,{\bf
c}})(v)$$ for all $i\in I$ and all $v\in V$. Since ${\mathcal
I}_H\subseteq \ker\,\rho$, all operators $\tilde{\rho}(u_i)$
annihilate the subspace $1\otimes V$ of $\widetilde{V}$. Since
$\widetilde{V}$ is generated by $1\otimes V$ as a $\g$-module and
the $\mathbb C$-span of the $u_i$'s is invariant under the adjoint
action of $\g$ on $U(\g)$, it must be that
$u_i\in\ker\tilde{\rho}$ for all $i\in I$. Since the $u_i$'s
generate the ideal $\mathcal I$, we obtain ${\mathcal I}\subseteq
\tilde{\mathcal I}$. Since the primitive ideals $\mathcal I$ and
$\tilde{\mathcal I}$ share the same associated variety,
[\cite{BK}, Corollar~3.6] implies that $\tilde{\mathcal
I}={\mathcal I}$, finishing the proof.
\end{proof}
By the main result of [\cite{DDeC}] the algebra $H_\chi$ is
isomorphic to the finite $W$-algebra $W^{\rm fin}(\g,e)$. Recall
that $W^{\rm fin}(\g,e)$ is a BRST quantisation of the Poisson
algebra $\gr\, H_\chi$, and it is also isomorphic to the Zhu
algebra of the vertex (affine) $W$-algebra associated with $\g$
and $e$; see [\cite{DeSK}, \cite{DDeC}, \cite{GG}, \cite{P02}] for
more detail. Vertex $W$-algebras and their Zhu algebras arise in
quantum field theory and are studied extensively in mathematical
physics.
\begin{corollary}\label{Walg}
All finite $W$-algebras $W^{\rm fin}(\g,e)$ possess finite
dimensional representations.
\end{corollary}
\begin{proof}
(1) Let ${\mathfrak X}_\lambda$ denote the set of all primitive
ideals of $U(\g)$ with infinitesimal character $\lambda\in\h^*$.
In order to prove the corollary we need to find a {\it rational}
infinitesimal character $\mu\in\h^*$ and a primitive ideal
${\mathcal I}\in{\mathfrak X}_\mu$ such that
$\mathcal{VA}({\mathcal I})=\overline{\mathcal O}_\chi$. Since
$W^{\rm fin}(\g,e)\cong H_\chi$ as algebras, our statement then
will follow from Theorem~\ref{main}. Let $W=\la
s_\alpha\,|\,\,\alpha\in\Phi^+\ra$ be the Weyl group of $\g$.
Given $\lambda\in\h^*$ we set
$\Phi_\lambda^+:=\{\alpha\in\Phi^+\,|\,\,\lambda(h_\alpha)\in\Z\}$
and denote by $W_\lambda$ the subgroup of $W$ consisting of all
$w\in W$ such that $w\lambda-\lambda$ is a sum of roots. The
subgroup $W_\lambda$ is often referred to as the {\it integral
Weyl group} of $\lambda$; it is generated by all reflections
$s_\alpha$ with $\alpha\in\Phi_\lambda^+$; see [\cite{Ja}, 2.5]
for example. Recall that $\lambda\in\h^*$ is called {\it regular}
if $\lambda(h_\alpha)\ne 0$ for all $\alpha\in\Phi$.

We shall rely on well-known results of Barbasch--Vogan
[\cite{BV1}, \cite{BV2}]. It follows from [\cite{BV2}, Thm~4.5 and
Cor.~4.4] that there exists $\mu'\in\h^*$ such that $\mu'+\rho$
regular and $\mathcal{VA}\big({\rm
Ann}_{U(\g)}\,L(\mu')\big)\,=\,\overline{\mathcal O}_{\chi}$; see
{\it loc.\!\,\,cit.} and [\cite{Ja1}, 9.12]. In order to finish
the proof we need to find a {\it rational} $\mu''\in \h^*$ such
that $\mu''+\rho$ is regular and $W_{\mu''}=\,W_{\mu'}$. Then
$\mathcal{VA}\big({\rm
Ann}_{U(\g)}\,L(w(\mu''+\rho)-\rho)\big)\,=\,\overline{\mathcal
O}_{\chi}$ for some $w\in W$, because the collections of
associated varieties attached to $\mathfrak{X}_{\mu'}$ and
$\mathfrak{X}_{\mu''}$ coincide setwise; see [BV2, Thm~1.1]. It is
probably known that such a $\mu''$ exists, but we could not find a
reference and shall include a brief proof for completeness:

\smallskip

\noindent (2) Let $\h_{\mathbb Q}^*$ be  the $\mathbb Q$-span of
$\Pi$ in $\h^*$. This is a vector space over $\mathbb Q$ with
basis $\varpi_1,\ldots,\varpi_\ell$, where
$\varpi_i(h_{\alpha_j})=\delta_{ij}$ for $1\le i,j\le \ell$. Since
there exist $\mu\in\h_{\mathbb Q}^*$ such that
$\mu(h_\alpha)\not\in\Z$ for all $\alpha\in\Phi^+$, we may assume
that $W_{\mu'}\ne \{1\}$. Then
$\Phi_{\mu'}^+=\{\beta_1,\ldots,\beta_{M}\}$ is a nonempty set.
For $1\le i\le M$ set $r_i:=\mu'(h_{\beta_i})$. All $r_i$'s are
integers and the non-homogeneous system of linear equations
\begin{eqnarray}\label{sys}
x_1\varpi_1(h_{\beta_1})+\cdots+x_\ell\varpi_\ell(h_{\beta_1})&=&r_1\nonumber\\
\quad\ \ \vdots \quad\quad\qquad\quad\qquad\qquad\quad\vdots\quad &&\,\vdots\\
 x_1\varpi_1(h_{\beta_{M}})+\cdots+x_\ell\varpi_\ell(h_{\beta_{M}})&=&r_{M}\nonumber
\end{eqnarray}
is consistent. Since all coefficients of the system are rationals,
$\h^*_{\mathbb Q}$ contains a basis of the solution space to the
corresponding homogeneous system, say $\mu_1,\ldots,\mu_{q}$. As
$r_i\in\Z$ for all $i$, we can also find $\mu_0\in\h_{\mathbb
Q}^*$ such that any solution to (\ref{sys}) has the form
$\mu_0+y_1\mu_1+\cdots+y_q\mu_q$ for some $y_1,\ldots,
y_q\in\mathbb C$.

Let $E_{\mathbb Q}$ and $E_{\mathbb R}$ denote the spans of
$\mu_1,\ldots,\mu_q$ over $\mathbb Q$ and $\mathbb R$,
respectively. If $E_{\mathbb Q}=0$, then $\mu'=\mu_0\in\h_{\mathbb
Q}^*$ and we are done. So assume $E_{\mathbb Q}\ne 0$. Write
$\Phi^+\setminus\Phi_{\mu'}^+=\{\beta_{N-M+1},\ldots,\beta_{N}\}$.
For $N-M< i\le N$ define $v_i\in (E_{\mathbb R})^*$ by setting
$v_i(\nu)=\,\nu(h_{\beta_i})$ for all $\nu\in E_{\mathbb R}$. For
$N-M< i\le N$ and $a\in\mathbb Q$ set $H_i(a):=\{\nu\in E_{\mathbb
R}\,|\,\,v_i(\nu)=a\}$, and define
$$\Omega\,:=\,\textstyle{\bigcup}_{\,i=N-M+1}^N\,\textstyle{\bigcup}_{\,n\in\Z}\,\,\,
H_{i}(-\mu_0(h_{\beta_i})+n).$$ If $v_i=0$ for some $i>N-M$, then
$(\mu_0+\nu)(h_{\beta_i})=\mu'(h_{\beta_i})\not\in\Z$ for all
$\nu\in E_{\mathbb R}$, implying that
$H_{i}(-\mu_0(h_{\beta_i})+n)=\emptyset$ for all $n\in\Z$. Hence
each nonempty $H_{i}(-\mu_0(h_{\beta_i})+n)$ is an affine
hyperplane in the Euclidean space $E_{\mathbb R}$. From this it
follows that the set $E_{\mathbb R}\setminus \Omega$ is nonempty
and open in the Euclidean topology of $E_{\mathbb R}$. As the
subset $E_{\mathbb Q}$ of $E_{\mathbb R}$ is dense in the
Euclidean topology of $E_{\mathbb R}$, there exists $\nu_0\in
E_{\mathbb Q}$ such that $v_i(\nu_0)\not\in-\mu_0(h_{\beta_i})+\Z$
for $N-M <i\le N$.

Now put $\mu'':=\mu_0+\nu_0$. Then $\mu''\in\h_{\mathbb Q}^*$, by
our earlier remarks in the proof, and $\mu''(h_{\beta_i})\not\in\Z$
for all $i>N-M$. Because $\mu''$ is a solution to (\ref{sys}), we
also have that $\mu''(h_{\beta_i})=\mu'(h_{\beta_i})$ for $1\le i\le
N-M$. This yields $\Phi_{\mu'}^+=\Phi_{\mu''}^+$ forcing
$W_{\mu'}=W_{\mu''}$. Since $\mu'+\rho$ is regular, so is
$\mu''+\rho$. This completes the proof.
\end{proof}

\end{document}